\documentclass[12pt]{article}
\usepackage{epic,eepic,epsf,epsfig}
\usepackage{amsfonts,srcltx,mathrsfs}
\usepackage{amsmath}
\usepackage{amssymb}
\usepackage{amsbsy}
\usepackage{amsfonts}
\usepackage{graphicx}

\newtheorem{theorem}{Theorem}[section]%
\newtheorem{lemma}[theorem]{Lemma}%
\newtheorem{prop}[theorem]{Proposition}%

\def\a{\alpha} \def\b{\beta} \def\g{\gamma}
\def\d{\delta}
\def\r{\rho}    
 \def\ld{\lambda}  \def\z{\zeta}
\def\G{\Gamma}

\def\oa{{\overline \a}}

\def \lg {\langle} \def\rg{\rangle}
\def \l {\langle} \def\r{\rangle}

\def\mz{{\mathbb Z}}
\def\f{\noindent}

\def\S{\hbox{\rm S}}
\def\A{\hbox{\rm A}}
\def\TL{\hbox{\rm $\Gamma$L}}
\def\AGL{\hbox{\rm AGL}}
\def\ASL{\hbox{\rm ASL}}
\def\ATL{\hbox{\rm A$\Gamma$L}}
\def\ASIL{\hbox{\rm A$\Sigma$L}}

\def\Aut{\hbox{\rm Aut}}

\def\Cay{\hbox{\rm Cay}}

\def\PSL{\hbox{\rm PSL}}
\def\PGL{\hbox{\rm PGL}}
\def\PSU{\hbox{\rm PSU}}
\def\PSp{\hbox{\rm PSp}}

\def\Dip{\hbox{\rm Dip}}

\def\GD {\hbox{\rm GD}}

\def\@makefnmark{}

\newcommand{\qed}{\hfill \mbox{\raisebox{0.7ex}{\fbox{}}} \vspace{4truemm}}

\def\demo{\f {\bf Proof.}\hskip10pt}

\begin{document}

\title{Pentavalent symmetric graphs of order twice a prime power} 
\author{Yan-Quan Feng, Jin-Xin Zhou  \\ {\small\em  Mathematics, Beijing Jiaotong University,} {\small\em Beijing 100044, P.R. China}\\ [2ex]
Yan-Tao Li  \\ {\small\em  College of Arts and Science, Beijing Union University,} {\small\em Beijing 100191, P.R. China}}

\footnotetext{E-mail addresses: yqfeng$@$bjtu.edu.cn, jxzhou@bjtu.edu.cn,
yantao@ygi.edu.cn}

\date{}
\maketitle \vspace{-.8cm}

\begin{abstract}
A connected symmetric graph of prime valency is {\em basic} if its
automorphism group contains no nontrivial normal subgroup having
more than two orbits. Let $p$ be a prime and $n$ a positive integer.
In this paper, we investigate properties of connected pentavalent
symmetric graphs of order $2p^n$, and it is shown that a connected
pentavalent symmetric graph of order $2p^n$ is basic if and only if
it is either a graph of order $6$, $16$, $250$, or a graph of three
infinite families of Cayley graphs on generalized dihedral groups --
one family has order $2p$ with $p=5$ or $5 \mid (p-1)$,  one family
has order $2p^2$ with $5 \mid (p\pm 1)$, and the other family has
order $2p^4$. Furthermore, the automorphism groups of these basic
graphs are computed. Similar works on cubic and tetravalent
symmetric graphs of order $2p^n$ have been done.

It is shown that basic graphs of  connected pentavalent symmetric
graphs of order $2p^n$ are symmetric elementary abelian covers of
the dipole $\Dip_5$, and with covering techniques, uniqueness and
automorphism groups of these basic graphs are determined. Moreover,
symmetric $\mz_p^n$-covers of the dipole $\Dip_5$ are classified. As
a byproduct,  connected pentavalent symmetric graphs of order $2p^2$
are classified.

\medskip
\f {\bf Key Words:} Symmetric graph, Cayley graph, regular covering, normal cover.\\
\f {\bf 2000 Mathematics Subject Classification:} 05C25, 20B25.

\end{abstract}

\section{Introduction}

Let $G$ be a permutation group on a set $\Omega$ and $\alpha\in
\Omega$. Denote by $G_{\alpha}$ the stabilizer of $\alpha$ in $G$,
that is, the subgroup of $G$ fixing the point $\alpha$. We say that
$G$ is {\em semiregular} on $\Omega$ if $G_{\alpha}=1$ for every
$\alpha \in \Omega$ and {\em regular} if $G$ is transitive and
semiregular. We will use the symbol $\mz_n$, both for the cyclic
group of order $n$ and for the ring of integers modulo $n$ (and for
the field of order $n$ if $n$ is a prime). Denote by $\mz_n^*$ the
multiplicative group of units of $\mz_n$, by $D_n$ the dihedral
group of order $2n$, and by $A_n$ and $S_n$ the alternating group
and the symmetric group of degree $n$, respectively.

All graphs in this article are finite, connected and simple,
unless explicitly stated. For a graph $\G$, let $V(\G)$, $E(\G)$ and
$\Aut(\G)$ denote the vertex set, edge set and full automorphism
group of $\G$, respectively. An {\em $s$-arc} in a graph $\G$ is an
ordered $(s+1)$-tuple $(v_0, v_1, \cdots ,  v_s)$ of $s+1$ vertices
such that  $\{v_{i-1}, v_i\}\in E(\G)$ for $1\leq i\leq s$ and
$v_{i-1}\neq v_{i+1}$ for $1\leq i \leq s-1$, and a $1$-arc is also
called an {\em arc}. For a subgroup $G$ of $\Aut(\G)$ of a graph
$\G$, the graph $\G$ is said to be $(G, s)$-{\em arc-transitive} or
$(G, s)$-{\em regular} if $G$ acts transitively or regularly on the
set of $s$-arcs of $\G$, and $(G, s)$-{\em transitive} if $G$ acts
transitively on the set of $s$-arcs but not on the set of
$(s+1)$-arcs of $\G$. A graph $\G$ is said to be $s$-{\em
arc-transitive}, $s$-{\em regular} or $s$-{\em transitive} if it is
$(\Aut(\G),s)$-arc-transitive, $(\Aut(\G), s)$-regular or
$(\Aut(\G),s)$-transitive. In particular, $0$-arc-transitive means
{\em vertex-transitive}, and $1$-arc-transitive means {\em
arc-transitive} or {\em symmetric}.

Let $\G$ be a graph and $N\leq \Aut(\G)$. The {\em quotient graph}
$\G_N$ of $\G$ relative to $N$ is defined as the graph
with vertices the orbits of $N$ on $V(\G)$ and with two orbits
adjacent if there is an edge in $\G$ between those two orbits. The
theory of quotient graph is widely used to investigate symmetric
graphs. Let $\G$ be a symmetric graph and $N\lhd
\Aut(\G)$. If $\G$ and $\G_N$ have same valency, the graph $\G$
is said to be a {\em normal cover} of $\G_N$ and the graph $\G_N$ is
said to be a {\em normal quotient} of  $\G$. In this case, $N$ is
semiregular on $V(\G)$. There are two steps to study a symmetric graph $\G$ --- the first step is to investigate normal quotient graph
$\G_N$ for some normal subgroup $N$ of $\Aut(\G)$ and the second step is to reconstruct
the original graph $\G$ from the normal quotient $\G_N$ by using covering
techniques. This is usually done by taking the normal subgroup $N$
as large as possible and then the graph $\G$ is reduced to a `basic graph'. The situation seems to be somewhat more promising with $2$-arc-transitive
graphs, and the strategy for the structural analysis of these graphs, based on taking normal
quotients, was first laid out by Praeger (see \cite{P1,P2,P3}). The strategy works for
locally primitive graphs, that is, vertex-transitive graphs with vertex stabilizers acting primitively on the corresponding neighbors¡¯ sets (see \cite{P4,P5}).

As for the first step, let us define some notations. A graph $\G$ is called {\em basic} if $\G$ has
no proper normal quotient. Then a locally primitive graph is basic if and only if it has no nontrivial normal subgroup having more than two orbits. A graph is {\em quasiprimitive} if every nontrivial normal subgroup of its automorphism group is transitive, and is {\em biquasiprimitive} if it has a nontrivial normal subgroup with two orbits but no such subgroup with more than two orbits. Therefore for locally primitive graphs, basic graphs are equivalent to quasiprimitive or biquasiprimitive graphs, which have received most of the attention thus far.
In \cite{IP}, Ivanov and Praeger
completed the classification of quasiprimitive $2$-arc-transitive graphs of affine type, and Baddeley gave a detailed
description of quasiprimitive $2$-arc-transitive graphs of twisted wreath type~\cite{B}. A similar description of
2-arc-transitive graphs associated with Suzuki groups and Ree groups was obtained by Fang and
Praeger~\cite{FP,FP2}. Classifications of quasiprimitive $2$-arc-transitive graphs of odd order and prime power order have been completed by Li~\cite{Li,Li2,Li3}, and based on this approach, finite vertex-primitive $2$-arc-regular graphs have been classified \cite{FLW} and finite $2$-arc-transitive Cayley graphs of abelian groups have been determined \cite{LP}. Most recently, symmetric graphs of diameter $2$ admitting an affine-type quasiprimitive group were investigated by Amarra et al.~\cite{AGP}, and
an infinite family of biquasiprimitive $2$-arc-transitive
cubic graphs were constructed by Devillers et al.~\cite{DGLP}.

Based on the stabilizers of pentavalent symmetric
graphs given by Guo and Feng~\cite{GF}, in this paper we prove that
normal quotient graphs of connected pentavalent symmetric graph of
order twice a prime power can be $K_6$, $FQ_4$ (the folded hypercube of order $16$ ),
${\cal CD}_{p}$ $(p=5$ or $5\ |\ (p-1))$, ${\cal CGD}_{p^2}^1$
$(p=5$ or $5\ |\ (p-1))$, ${\cal CGD}_{p^2}^2$ $(5\ |\ (p\pm 1))$,
${\cal CGD}_{p^3}$ $(p=5$ or $5\ |\ (p-1))$, or ${\cal CGD}_{p^4}$,
where the graphs are defined in Eqs~(\ref{eq-0})-(\ref{eq-4}). Automorphism groups of these normal quotients are computed, and among them, basic ones are determined, which are $K_6$, $FQ_4$, ${\cal CGD}_{5^3}$, ${\cal CD}_{p}$ $(p=5$ or $5\ |\ (p-1))$, ${\cal CGD}_{p^2}^2$ $(5\ |\
(p\pm 1))$ and ${\cal CGD}_{p^4}$ $(p=3$ or $p\geq 7)$. Similar works
on cubic symmetric graphs and tetravalent $2$-arc-transitive graphs of order twice a prime power were done by Kwak and the first two authors~\cite{FK6,ZF}.

As for the second step, regular covering (for
notation, see Section~\ref{s3}) is becoming an active topic in
algebraic graph theory. In~\cite{DMW}, regular covers of complete graphs whose group of covering transformations is either cyclic or isomorphic to $\mz_p^2$, $p$ a prime, and whose fibre-preserving subgroup of automorphisms acts $2$-arc-transitively, were classified. This result has been extended to the case
where the group of covering transformations is isomorphic to $\mz_p^3$, p a prime \cite{DKX}.
Some general methods of elementary abelian
coverings were developed in~\cite{DKX1,MMP,MMP1}. By using the
method developed in~\cite{MMP1}, Malni\v{c} and
Poto\v{c}nik~\cite{MP0} classified all vertex-transitive elementary
abelian covers of the Petersen graph. Symmetric cyclic or elementary
abelian covers of the complete graph $K_4$, the complete bipartite
graph $K_{3,3}$, the cube $Q_3$ and the Petersen graph $O_3$, were
classified in \cite{FK3,FK4,FK5,FKW,FW}. Symmetric
elementary abelian covers of the unique connected cubic symmetric
graph of order $14$, $16$ or $18$ were classified  in~\cite{O1,O2,O3}. By
using the above covers, together with group theory techniques,  many
classifications of symmetric graphs have been obtained -- for
example, symmetric cubic graphs of order $rp$ or $tp^2$ were
classified for each $2\leq r\leq 16$ and $2\leq t\leq 10$.
Classification of symmetric graphs with a given order has been
widely investigated, and for more results, see \cite{Chao, CO,LLM,
PWX1,PWX2,WX}. In the above papers, graphs and their
covers are simple, that is, no loops and multiple edges. Regular
covers of non-simple graphs were also considered in literature
and in this case, automorphism groups of  non-simple graphs are
usually considered as permutation groups on the sets of arcs of these
graphs. For example, to classify tetravalent non-Cayley graph of
order four times a prime, Zhou~\cite{Zhou} considered
vertex-transitive covers of non-simple graphs of order $4$.

 To determine the uniqueness of normal quotient graphs of connected pentavalent symmetric graph of order twice a prime power for some given orders and to compute their automorphism groups, covering techniques are employed.  In this paper we first prove that these normal quotients are symmetric elementary abelian covers of the dipole $\Dip_5$ and then determine all symmetric elementary abelian covers of  $\Dip_5$, which consist of four infinite families of Cayley graphs on
generalized dihedral groups, that is, the graphs ${\cal
CGD}_{p^2}^1$  ($p=5$ or $5\ |\ (p-1)$), ${\cal CGD}_{p^2}^2$ ($5\
|\ (p\pm 1)$), ${\cal CGD}_{p^3}$ ($p=5$ or $5\ |\ (p-1)$) and
${\cal CGD}_{p^4}$. These covers are not isomorphic to each other and their full automorphism groups are computed.
As an application, pentavalent symmetric graphs of order
twice a prime square are classified.

\section{Preliminaries}

In this section, we describe some preliminary results which will be
used later. First we describe stabilizers of connected pentavalent
symmetric  graphs.

\begin{prop}\label{stabilizer-val5}{\rm\cite[Theorem~1.1]{GF}}
Let $X$ be a connected pentavalent $(G,s)$-transitive graph for some
$G\leq \Aut(X)$ and $s\geq 1$. Let $v\in V(X)$. Then $s\leq 5$ and
one of the following holds:

\begin{enumerate}\setlength{\parsep}{-1pt}\setlength{\itemsep}{-1pt}

\item[{\rm (1)}] For $s=1$, $G_v\cong\mz_5$, $D_5$ or $D_{10}$;

\item[{\rm (2)}] For $s=2$, $G_v\cong F_{20}$,
$F_{20}\times\mz_2$, $\A_5$ or $\S_5$, where $F_{20}$ is the
Frobenius group of order $20$;

\item[{\rm (3)}] For $s=3$, $G_v\cong F_{20}\times\mz_4$,
$\A_4\times\A_5$, $\S_4\times\S_5$ or $(\A_4\times\A_5)\rtimes\mz_2$
with $\A_4\rtimes\mz_2=\S_4$ and $\A_5\rtimes\mz_2=\S_5$;

\item[{\rm (4)}] For $s=4$, $G_v\cong\ASL(2,4)$, $\AGL(2,4)$,
$\ASIL(2,4)$ or $\ATL(2,4)$;

\item[{\rm (5)}] For $s=5$,
$G_v\cong\mz_2^6\rtimes\TL(2,4)$.

\end{enumerate}
\end{prop}

For a subgroup $H$ of a group $G$, denote by $C_G(H)$ the
centralizer of $H$ in $G$ and by $N_G(H)$ the normalizer of $H$ in $G$. Then $C_G(H)$ is normal in $N_G(H)$.

\begin{prop}{\rm\cite[Theorem~6.11]{Hup}} \label{NC}
The quotient group $N_G(H)/C_G(H)$ is isomorphic to a subgroup of
the automorphism group $\Aut(H)$ of $H$.
\end{prop}

Let $G$ be a finite group and let $\pi(G)=\{p\ |\ p \mbox{ is a
prime divisor of $|G|$}\}$. Herzog~\cite{Herzog} and Shi~\cite{Shi}
(also see \cite{HL}) classified nonabelian finite simple groups $G$
for $|\pi(G)|=3$ and $|\pi(G)|=4$ respectively, from which one may
deduce the following proposition.

\begin{prop} \label{non-abelian-simple-groups}
Let $p\geq 7$ be a prime and let $G$ be a nonabelian
simple group.

\begin{enumerate}
  \item [{\rm (1)}] If $|\pi(G)|=3$ then $G\cong A_5,A_6,\PSL(2,7),\PSL(2,8),\PSL(2,17),\PSL(3,3),
      \PSU(3,3)$ or $\PSU(4,2)$ with order $2^2\cdot 3\cdot 5$, $2^3\cdot 3^2\cdot 5$, $2^3\cdot 3\cdot 7$, $2^3\cdot 3^2\cdot 7$, $2^4\cdot 3^2\cdot 17$, $2^4\cdot 3^3\cdot 13$, $2^5\cdot 3^3\cdot 7$ or $2^6\cdot 3^4\cdot 5$, respectively.

\item [{\rm (2)}] Let $\pi(G)=\{2,3,5,p\}$ with $p\geq 7$. If $p^2\mid |G|$, $2^{11}\nmid |G|$ and $3^3\nmid |G|$, then  $G\cong \PSL(2,49)$ or $\PSp(4,7)$ with order $2^4\cdot 3\cdot 5^2\cdot 7^2$ or $2^8\cdot 3^2\cdot 5^2\cdot 7^4$, respectively.
\end{enumerate}
\end{prop}

Let $G$ be a finite group and $S$ a subset of $G$ with $1\not\in S$
and $S^{-1}=S$. The {\em Cayley graph} $\G=\Cay(G,S)$ on $G$ with
respect to $S$ is defined to have vertex set $V(\G)=G$ and edge set
$E(\G)=\{\{g,sg\}\ |\ g\in G,s\in S\}$. It is well-known that
$\Aut(\G)$ contains the right regular representation $R(G)$ of $G$,
the acting group of $G$ by right multiplication, and $\G$ is
connected if and only if $G=\langle S \rangle$, that is, $S$
generates $G$. A Cayley graph $\Cay(G,S)$ is said to be {\em normal}
if the right regular representation $R(G)$ of $G$ is normal in
$\Aut(\Cay(G,S))$. By Godsil~\cite{G} or Xu~\cite{X}, we have the
following result.

\begin{prop} \label{stabilizer}
Let $\G=\Cay(G,S)$ be a connected Cayley graph on a finite group $G$ with respect to $S$, and let $A=\Aut(\G)$. Then
$N_A(R(G))=R(G)\rtimes\Aut(G,S)$, where $N_A(R(G))$ is the normalizer of $R(G)$ in $A$. In particular, $\G$ is normal if and only if $A_1={\Aut}(G,S)$.
\end{prop}

Let $p$ be a prime and let $D_p=\langle a, b \ | \ a^{p}=b^{2}=1,
b^{-1}ab=a^{-1}\rangle $ be the dihedral group of order $2p$. For
$p=5$, let $\ell=1$ and for $5\ |\ (p-1)$, let $\ell$ be an element
of order $5$ in $\mz_p^*$. Define

\begin{equation}\label{eq-0}
{\cal CD}_{p}=\Cay(D_{p}, \{b, ab, a^{\ell+1}b,
a^{\ell^{2}+\ell+1}b, a^{\ell^{3}+\ell^{2}+\ell+1}b\}).
\end{equation}

\f The graph ${\cal CD}_{p}$ is symmetric because the map $\a: a
\mapsto a^{\ell}, b \mapsto ab$ induces an automorphism of $D_p$
permuting the elements in $\{b, ab, a^{\ell+1}b, a^{\ell^{2}+\ell+1}b,
a^{\ell^{3}+\ell^{2}+\ell+1}b\}$ cyclicly. By
\cite[Theorem~3.1]{FL}, ${\cal CD}_{p}$ is independent of the choice
of $\ell$, and by \cite{CO} and \cite[Theorem~3.1]{FL}, we have the
following.

\begin{prop}\label{2p}
Let $\G$ be a connected  pentavalent  edge-transitive graph of order $2p$ for a prime $p$. Then $\G$ is symmetric and one of the following holds:
\begin{enumerate}

\item [{\rm (1)}] $X\cong K_6$, the complete graph of order $6$ and $\Aut(K_6)\cong S_6$;
\item [{\rm (2)}] $X\cong {\cal CD}_{5}$ $( \cong K_{5,5})$, the complete bipartite graph of order $10$ and $\Aut({\cal CD}_{5})\cong (S_5\times S_5)\rtimes\mz_2$;
\item [{\rm (3)}]  $X\cong {\cal CD}_{p}$ with $5\ |\ (p-1)$. For $p=11$, $\Aut({\cal
CD}_{p})\cong \PGL(2, 11)$ and for $p\geq 31$, $\Aut({\cal
CD}_{p})\cong D_p\rtimes\mz_5$.
\end{enumerate}
\end{prop}

Let $p$ be a prime and let $D_{p^2}=\langle a, b \ | \
a^{p^2}=b^{2}=1, b^{-1}ab=a^{-1}\rangle $ be the dihedral group of
order $2p^2$. Let $5\ |\ (p-1)$ and $\ell$ an element of order $5$
in $\mz_{p^2}^*$. Define

\begin{equation}\label{eq-00}
{\cal CD}_{p^2}=\Cay(D_{p^2}, \{b, ab, a^{\ell+1}b,
a^{\ell^{2}+\ell+1}b, a^{\ell^{3}+\ell^{2}+\ell+1}b\}).
\end{equation}

\f Similar to the graph ${\cal CD}_{p}$, the graph ${\cal CD}_{p^2}$
is symmetric because the map $a \mapsto a^{\ell}, b \mapsto ab$
induces an automorphism of $D_{p^2}$. By \cite[Theorem~3.1]{FL},
${\cal CD}_{p^2}$ is independent of the choice of $\ell$, and by
\cite[Proposition~2.2 and Theorem~A]{KKO}, we have the following.

\begin{prop}\label{normal2p2}
Let $p$ be a prime and $\G$ be a connected pentavalent Cayley graph
on $D_{p^2}$. If $\G$ is $N_{\Aut(\G)}(R(D_{p^2}))$-arc-transitive,
then $5\ |\ (p-1)$ and $\G\cong {\cal CD}_{p^2}$ with $\Aut({\cal
CD}_{p^2})\cong R(D_{p^2}) \rtimes \mz_5$.
\end{prop}

Let $\G$ be a connected symmetric graph of prime valency and let
$G\leq \Aut(\G)$ be arc-transitive. Let $N$ be a normal subgroup of
$G$. In view of \cite[Theorem 9]{Lor}, we have:

\begin{prop}\label{QG}
If $N$ has more than two orbits then the quotient graph $\G_{N}$ has
the same valency as $\G$ and $N$ is the kernel of $G$ acting on the
set of orbits of $N$. Furthermore, $N$ is semiregular on $V(\G)$ and
$\G_{N}$ is $G/N$-arc-transitive.
\end{prop}

For an abelian group $H$, the generalized dihedral group of $H$,
denoted by $\GD_H$, is the semidirect product  $H \rtimes \mz_2$
with the involution of $\mz_2$ inverting every element in $H$. To
end this section, we consider some special bipartite graphs.

\begin{lemma}\label{semiregulargraph}
Let $\G$ be a bipartite graph and $H$ an abelian semiregular
automorphism group of $\G$ with the two bipartite sets of $\G$ as its
orbits. Then $\G$ is a Cayley graph on $\GD_H$.
\end{lemma}

\demo Let $B_1$ and $B_2$ be the bipartite sets of $\G$. Then $H$
acts regularly on each of $B_1$ and $B_2$, and we may assume that
$B_1=\{h\ |\ h\in H\}$ and $B_2=\{h'\ |\ h\in H\}$. The actions of
$H$ on $B_1$ and $B_2$ are just by right multiplication, that is,
$h^g=hg$ and $(h')^g=(hg)'$ for any $h, g\in H$. Let the neighbors
of $1$ in $\G$ be $h_1'$, $h_2'$, $\cdots$, $h_n'$, where
$h_1,h_2,\cdots,h_n\in H$. Since $H$ is abelian, for any $h\in H$,
the neighbors of $h$ are $(hh_1)'$, $(hh_2)'$, $\cdots$, $(hh_n)'$,
and furthermore, the neighbors of $h'$ are $hh_1^{-1}$, $hh_2^{-1}$,
$\cdots$, $hh_n^{-1}$. It is easy to check that the map $\a$,
defined by $h\mapsto (h^{-1})'$, $h'\mapsto h^{-1}$, $h\in H$, is an
automorphism of $\G$ of order $2$. Now for any $g,h\in H$, we have
$h^{g\a}=(g^{-1}h^{-1})'=h^{\a g^{-1}}$ and
$(h')^{g\a}=g^{-1}h^{-1}=(h')^{\a g^{-1}}$. It follows that $g\a=\a
g^{-1}$, that is, $\a^{-1} g\a=g^{-1}$. Thus, $\langle
H,\a\rangle=GD_H$ and $\G$ is a Cayley graph on $GD_H$. \hfill\qed

\section{Cayley graphs on $\GD_{\mz_p^n}$ with $n=2,3,4$\label{s3}}

In this section, we shall construct four families of Cayley graphs
on generalized dihedral groups, and investigate their automorphisms,
which will be used later.
For a prime $p$, write $\GD_{p^2}=\GD_{\mz_{p}^{2}}$,
$\GD_{p^3}=\GD_{\mz_{p}^{3}}$ and $\GD_{p^4}=\GD_{\mz_{p}^{4}}$ for
short.

Let $\GD_{p^{2}}=\l a, d, h \ | \ a^p=d^p=h^2=[a,d]=1,
h^{-1}ah=a^{-1}, h^{-1}dh=d^{-1}\r$. For $p=5$, let $\ell=1$, and
for $5\ |\ (p-1)$, let $\ell$ be an element of order $5$ in
$\mz_p^*$. Define
\begin{equation}
{\cal CGD}_{p^2}^1=\Cay(\GD_{p^2}, \{h, ah,
a^{\ell(\ell+1)^{-1}}d^{\ell^{-1}}h, a^{\ell}d^{(\ell+1)^{-1}}h, dh
\}).\label{eq-1}
\end{equation}

For $5\ |\ (p\pm 1)$, let $\ld$ be an element in $\mz_p^*$ such that
$\ld^2=5$. Define
\begin{equation}
{\cal CGD}_{p^2}^2=\Cay(\GD_{p^2}, \{h, ah, a^{2^{-1}(1+\ld)}dh,
ad^{2^{-1}(1+\ld)}h, dh\}).\label{eq-2}
\end{equation}

Let $\GD_{p^{3}}=\l a , b, d, h \ | \
a^p=b^p=d^p=h^2=[a,b]=[a,d]=[b,d]=1, h^{-1}ah=a^{-1},
h^{-1}bh=b^{-1}, h^{-1}dh=d^{-1} \r$. For $p=5$, let $\ell=1$, and
for $5\ |\ (p-1)$, let $\ell$ be an element of order $5$ in
$\mz_p^*$. Define
\begin{equation}
{\cal CGD}_{p^3}=\Cay(\GD_{p^3}, \{h, ah, bh,
a^{-\ell^2}b^{-\ell}d^{-\ell^{-1}}h, dh\}).\label{eq-3}
\end{equation}

Let $\GD_{p^{4}}=\l a , b, c, d, h \ | \ a^p=b^p=c^p=d^p=h^2=[a,
b]=[a, c]=[a, d]=[b, c]=[b, d]=[c, d]=1, h^{-1}ah=a^{-1},
h^{-1}bh=b^{-1}, h^{-1}ch=c^{-1}, h^{-1}dh=d^{-1} \r$. Define
\begin{equation}
{\cal CGD}_{p^4}=\Cay(\GD_{p^4},\{h, ah, bh, ch,dh\}).\label{eq-4}
\end{equation}

\begin{theorem} \label{lem-3.1}
Let $\G=\Cay(G,S)$ be one of the graphs defined in
Eqs~(\ref{eq-1})-(\ref{eq-4}). Let $P$ be a Sylow $p$-subgroup of
$R(G)$ and let $A=\Aut(\G)$. Then $\G$ is $N_A(P)$-arc-transitive.

\begin{enumerate}
\item [{\rm (1)}] Let $\G={\cal CGD}_{p^2}^1$ $(p=5$ or $5\mid (p-1))$. If $p=5$
then $N_A(R(GD_{5^2}))\cong R(GD_{5^2})\rtimes F_{20}$ and if $5\ |\
(p-1)$ then $N_A(R(GD_{p^2}))\cong R(GD_{p^2})\rtimes\mz_5$.
Furthermore, $|N_A(P)|\not=20p^2$.

  \item [{\rm (2)}] Let $\G={\cal CGD}_{p^2}^2$ $(5\ |\
  (p\pm 1))$. Then $N_A(R(GD_{p^2}))\cong R(GD_{p^2})\rtimes D_5$ and $|N_A(P)|$ has a divisor $20p^2$.

 \item [{\rm (3)}] Let $\G={\cal CGD}_{p^3}$ $(p=5$ or $5\mid (p-1))$. If $p=5$
  then $N_A(R(GD_{5^3}))\cong R(GD_{5^3})\rtimes S_5$ and if $5\ |\ (p-1)$ then $R(GD_{p^3})\rtimes\mz_5\leq N_A(R(GD_{p^3}))$.

 \item [{\rm (4)}] Let $\G={\cal CGD}_{p^4}$. Then  $N_A(R(GD_{p^4}))\cong R(GD_{p^4})\rtimes S_5$.

\end{enumerate}
\end{theorem}

\demo  By Proposition~\ref{stabilizer}, $N_A(R(G))=R(G)\rtimes
\Aut(G,S)$. For each graph in
Eqs~(\ref{eq-1})-(\ref{eq-4}), we have $\langle S\rangle =G$ and
$|R(G):P|=2$. Thus, $\G$ is connected and $\Aut(G,S)$ acts
faithfully on $S$, implying $\Aut(G,S)\leq S_5$. Furthermore, $P\lhd
R(G)$ and hence $P\lhd N_A(R(G))$, forcing that $N_A(R(G))\leq
N_A(P)$. If (1)-(4) are true, then $\G$ is $N_A(P)$-arc-transitive because it is $N_A(R(G))$-arc-transitive. To finish the proof, it suffices to show (1)-(4).

Let $\G={\cal CGD}_{p^2}^1=\Cay(G, S)$. Then $p=5$ or $5\mid (p-1)$.
Furthermore, $G=\GD_{p^{2}}=\l a, d, h \ | \ a^p=d^p=h^2=[a,d]=1,
h^{-1}ah=a^{-1}, h^{-1}dh=d^{-1}\r$ and $S=\{h, ah,
a^{\ell(\ell+1)^{-1}}d^{\ell^{-1}}h, a^{\ell}d^{(\ell+1)^{-1}}h, dh
\}$, where $\ell=1$ for $p=5$ and $\ell$ is an element of order $5$
in $\mz_p^*$ for $5\ |\ (p-1)$. In particular, $\ell^5=1$,
$\ell^4+\ell^3+\ell^2+\ell+1=0$ and $(\ell+1)^{-1}=-\ell^3-\ell$ in
the field $\mz_p$. It is easy to check that $\a: h \mapsto ah, a
\mapsto a^{\ell(\ell+1)^{-1}-1}d^{\ell^{-1}}, d \mapsto a^{-1}$
induces an automorphism of order $5$ of $GD_{p^2}$ permuting
the elements in $S$ cyclicly.

For $p=5$, the map $\b: a\mapsto ad^3$, $d\mapsto a^3d$ and
$h\mapsto h$, induces an automorphism of order $4$ of $GD_{p^2}$
permuting the elements in $\{ah,ad^3h,dh,a^3dh\}$ cyclicly. Take $\g\in
\Aut(G,S)$ such that $h^\g=h$ and $(ah)^\g=ah$. Then
$\{a^3d,ad^3,d\}^\g=\{a^3d,ad^3,d\}$. If $d^\g=a^3d$ or $ad^3$ then
$(a^3d)^d=ad$ or $a^4d^3$, which is impossible because
$ad,a^4d^3\not\in \{a^3d,ad^3,d\}$. Thus, $d^\g=d$ and hence $\g=1$,
implying $|\Aut(G,S)|=20$. By Lemma~\ref{stabilizer-val5},
$N_A(R(GD_{5^2}))\cong R(GD_{5^2})\rtimes F_{20}$. Since
$N_A(R(G))\leq N_A(P)$, we have $|N_A(P)|\not=20p^2$.

For $5\ |\ (p-1)$, take $\d\in \Aut(G,S)$ such that $h^\d=h$. Then
$\d$ fixes $\{a, a^{\ell(\ell+1)^{-1}}d^{\ell^{-1}}$,
$a^{\ell}d^{(\ell+1)^{-1}}$, $d\}$ setwise, and hence fixes $a\cdot
a^{\ell(\ell+1)^{-1}}d^{\ell^{-1}}\cdot
a^{\ell}d^{(\ell+1)^{-1}}\cdot d$, that is,
$a^{\ell^3+2\ell+2}d^{2\ell^4+\ell^2+2}$. Thus,
$(a^rd^{r\ell^4})^\d=a^rd^{r\ell^4}$, where $r=\ell^3+2\ell+2$.
Suppose $r=0$. Then $2\ell^3+2\ell^2+1=\ell^2r=0$ and
$3\ell^3+4\ell^2-2\ell=2(2\ell^3+2\ell^2+1)-r=0$, that is,
$3\ell^2+4\ell-2=0$. On the other hand,
$2\ell^4+\ell^3+2\ell^2=2(\ell^4+\ell^3+\ell^2+\ell+1)-r=0$, that
is, $2\ell^2+\ell+2=0$. It follows that
$5\ell(\ell+1)=(3\ell^2+4\ell-2)+(2\ell^2+\ell+2)=0$, and hence
$p=5$, a contradiction. This implies that
$(ad^{\ell^4})^\d=ad^{\ell^4}$. One may compute the following
equations.

\parbox{6cm}{
\begin{eqnarray*}
& a^{ (1-\ell(1+\ell)^{-1})}=(ad^{\ell^4})\cdot
(a^{\ell(\ell+1)^{-1}}d^{\ell^{-1}})^{-1}, &
a^{-1}=(ad^{\ell^4})\cdot (a^{\ell}d^{(\ell+1)^{-1}})^{-
\ell^4(\ell+1)},\\&
a^{-\ell(\ell+1)^{-1}}=(ad^{\ell^4})(a^{1+\ell(\ell+1)^{-1}}d^{\ell^{-1}})^{-1},
&a^{(1-\ell^4)}=(ad^{\ell^4})(ad)^{-\ell^4},\\
&a^{(\ell^3+\ell^2)}=(ad^{\ell^4})(a^{1+\ell}d^{(\ell+1)^{-1}})^{-\ell^4(\ell+1)},&a=(ad^{\ell^4})d^{-\ell^4}.
\end{eqnarray*}}\hfill
\parbox{1.0cm}{
\begin{eqnarray*}
 \\
 \hfill (*)\\
 \end{eqnarray*}}

Recall that $\d$ fixes $S_1=\{a, a^{\ell(\ell+1)^{-1}}d^{\ell^{-1}},
a^{\ell}d^{(\ell+1)^{-1}}, d\}$ setwise. If $\d$ fixes $a$ then it
also fixes $d$ because $d=(ad^{\ell^4})^\ell a^{-\ell}$, which
implies that $\d=1$. By Eq($*$), if $\d$ fixes any element in $S_1$
then $\d$ fixes $a$ and therefore $\d=1$.

Suppose that $a^\d\not=a$. Then either $\langle \d\rangle$ has two
orbits of length $2$, or it is transitive on $S_1$. For the former,
$\d$ fixes one element in $\{a\cdot
a^{\ell(\ell+1)^{-1}}d^{\ell^{-1}}, a\cdot
a^{\ell}d^{(\ell+1)^{-1}}, a\cdot d\}$, that is, $\{
a^{1+\ell(\ell+1)^{-1}}d^{\ell^{-1}},  a^{1+\ell}d^{(\ell+1)^{-1}},
ad \}$. By Eq($*$), $\d$ fixes $a$, a contradiction. For the latter,
$\d^2$ has two orbits of length $2$ on $S_1$, and hence $\d^2$ fixes
one element in $\{a\cdot a^{\ell(\ell+1)^{-1}}d^{\ell^{-1}}, a\cdot
a^{\ell}d^{(\ell+1)^{-1}}, a\cdot d\}=\{
a^{1+\ell(\ell+1)^{-1}}d^{\ell^{-1}},  a^{1+\ell}d^{(\ell+1)^{-1}},
ad \}$. Again by Eq($*$), $\d^2$ fixes $a$ and $\d^2=1$, which is
impossible because $\d$ is transitive on $S_1$.

Since $\d=1$, we have that $|\Aut(G,S)|=5$, $N_A(R(GD_{p^2}))\cong
R(GD_{p^2})\rtimes\mz_5$ and $|N_A(R(GD_{p^2}))|=10p^2$. Suppose
$|N_A(P)|=20p^2$. Since $N_A(R(GD_{p^2}))\leq N_A(P)$, we have
$|N_A(P):N_A(R(GD_{p^2}))|=2$, implying $N_A(R(GD_{p^2}))\lhd
N_A(P)$. Note that $R(GD_{p^2})$ is a hall $\{2,p\}$-subgroup of
$N_A(R(GD_{p^2}))$. Then $R(GD_{p^2})$ is characteristic in
$N_A(R(GD_{p^2}))$ and hence $R(GD_{p^2})\lhd N_A(P)$. It follows
that $N_A(P)\leq N_A(R(GD_{p^2}))$, which is impossible. Thus,
$|N_A(P)|\not=20p^2$, completing the proof of (1).

Let $\G={\cal CGD}_{p^2}^2=\Cay(G, S)$. Then $5\ |\
  (p\pm 1)$, $G=GD_{p^2}=\l
a, d, h \ | \ a^p=d^p=h^2=[a,d]=1, h^{-1}ah=a^{-1},
h^{-1}dh=d^{-1}\r$ and $S=\{h, ah, a^{2^{-1}(1+\ld)}dh,
ad^{2^{-1}(1+\ld)}h, dh\}$, where $\ld^2=5$ in the field $\mz_p$.
The map $\a: h \mapsto ah, a \mapsto a^{2^{-1}(1+\ld)-1}d, d \mapsto
a^{-1}$ induces an automorphism of $\GD_{p^2}$ permuting the
elements in $S$ cyclicly. Take $\b\in \Aut(G,S)$ such that $h^\b=h$.
Then $\b$ fixes $\{a, a^{2^{-1}(1+\ld)}d, ad^{2^{-1}(1+\ld)}, d\}$
setwise and hence $a\cdot a^{2^{-1}(1+\ld)}d\cdot
ad^{2^{-1}(1+\ld)}\cdot d$, that is, $(ad)^{2+2^{-1}(1+\ld)}$. If
$2+2^{-1}(1+\ld)=0$ then $\ld=-5$, which is impossible because
$\ld^2=5$. Thus, $\b$ fixes $ad$.

Note that $\b$ fixes $\{a, a^{2^{-1}(1+\ld)}d, ad^{2^{-1}(1+\ld)},
d\}$  setwise. Clearly, if $\b$ fixes $a$ then $\b=1$. Assume that
$a^\b\not=a$. Then  $a^\b=a^{2^{-1}(1+\ld)}d, ad^{2^{-1}(1+\ld)}$ or
$d$. If $a^\b=a^{2^{-1}(1+\ld)}d$ then
$ad=(ad)^\b=a^{2^{-1}(1+\ld)}dd^\b$, that is, $d^\b=
a^{1-2^{-1}(1+\ld)}$. This is impossible because $d^\b\not\in \{a,
a^{2^{-1}(1+\ld)}d, ad^{2^{-1}(1+\ld)}, d\}$. Thus,
$a^\b\not=a^{2^{-1}(1+\ld)}d$, and similarly,
$a^\b\not=ad^{2^{-1}(1+\ld)}$ because otherwise
$d^\b=d^{1-2^{-1}(1+\ld)}$. It follows that $a^\b=d$ and
$d^\b=(a^{-1}ad)^\b=d^{-1}ad=a$, that is, $\b$ is an automorphism of
order $2$ of $GD_{p^2}$ induced by $a\mapsto d$, $d\mapsto a$ and
$h\mapsto h$. This implies that the subgroup of $\Aut(G,S)$ fixing
$h$ is $\langle \b\rangle$. Thus, $|\Aut(G,S)|=10$ and by
Proposition~\ref{stabilizer-val5}, $N_A(R(GD_{p^2}))\cong
R(GD_{p^2})\rtimes D_5$. Since $N_A(R(G))\leq N_A(P)$, we have
 that $|N_A(P)|$ has a divisor $20p^2$. This completes the proof of (2).

Let $\G={\cal CGD}_{p^3}=\Cay(G, S)$. Then $p=5$ or $5\mid (p-1)$.
Furthermore, $G=\GD_{p^{3}}=\l a , b, d, h \ | \
a^p=b^p=d^p=h^2=[a,b]=[a,d]=[b,d]=1, h^{-1}ah=a^{-1},
h^{-1}bh=b^{-1}, h^{-1}dh=d^{-1} \r$ and $S=\{h, ah, bh,
a^{-\ell^2}b^{-\ell}d^{-\ell^{-1}}h, dh\}$, where $\ell=1$ for $p=5$
and $\ell$ is an element of order $5$ in $\mz_p^*$ for $5\ |\ (p-1)$.
The map $\a: h \mapsto ah, a \mapsto ba^{-1}, b \mapsto
a^{-\ell^2-1}b^{-\ell}d^{-\ell^{-1}}, d \mapsto a^{-1}$ induces an
automorphism of $GD_{p^3}$ permuting the elements in $S$
cyclicly. Thus, $R(GD_{p^3})\rtimes\mz_5\leq N_A(R(GD_{p^3}))$. For
$p=5$, the map $\b: a\mapsto b$, $b\mapsto d$, $d\mapsto
a^{-1}b^{-1}d^{-1}$ and $h\mapsto h$, induces an automorphism of
order $4$ of $GD_{5^3}$ permuting the elements in
$\{ah,bh,dh,a^{-1}b^{-1}d^{-1}h\}$ cyclicly. Furthermore, any
permutation on $\{a,b,d\}$ with $h\mapsto h$ induces an automorphism
of $GD_{5^3}$. Thus, $\Aut(G,S)\cong S_5$ and $N_A(R(GD_{5^3}))\cong
R(GD_{5^3})\rtimes S_5$. This completes the proof of (3).

Let $\G={\cal CGD}_{p^4}=\Cay(G, S)$. Then $S=\{h, ah, bh, ch,dh\}$
and  $G=\GD_{p^4}=\l a , b, c, d, h \ | \ a^p=b^p=c^p=d^p=h^2=[a,
b]=[a, c]=[a, d]=[b, c]=[b, d]=[c, d]=1, h^{-1}ah=a^{-1},
h^{-1}bh=b^{-1}, h^{-1}ch=c^{-1}, h^{-1}dh=d^{-1} \r$. The map $\a:
h \mapsto ah, a \mapsto ba^{-1}, b \mapsto ca^{-1}, c \mapsto
da^{-1}, d \mapsto a^{-1}$ induces an automorphism of $GD_{p^3}$ permuting the elements in $S$ cyclicly. Furthermore, any
permutation on $\{a,b,c,d\}$ with $h\mapsto h$ induces an
automorphism of $GD_{p^4}$. Thus, $\Aut(G,S)\cong S_5$ and
$N_A(R(GD_{p^4}))\cong R(GD_{p^4})\rtimes S_5$. This completes the
proof of (4). \hfill\qed

\section{Symmetric elementary abelian covers of $\Dip_5$\label{s3}}

The main purpose of this paper is to determine basic graphs of pentavalent symmetric graphs of order twice a prime power and to compute their automorphism groups. To do that, we need  the so
called regular covering and this was treated in an extensive way~\cite{M,MMP,MNS,MP}. In this  section, we classify symmetric elementary abelian covers of the dipole $\Dip_5$ and compute their automorphism groups by combining Theorem~\ref{lem-3.1}.

An epimorphism ${\cal
P}:\widetilde{\Gamma} \mapsto \Gamma$ of graphs is called a {\em
regular covering projection} or {\em regular $N$-covering
projection} if $\Aut(\widetilde{\Gamma})$ has a semiregular subgroup
$N$ whose orbits on $V(\widetilde{\Gamma})$ coincide with the {\em
vertex fibres} ${\cal P}^{-1}(v), v\in V(\Gamma)$, and the arc and
edge orbits of $N$ coincide with the {\em arc fibres} ${\cal
P}^{-1}((u, v)), u \sim v$, and {\em edge fibres } ${\cal
P}^{-1}(\{u, v\}), u \sim v$, respectively. In particular, we call
the graph $\widetilde{\Gamma}$ a {\em regular cover} or an {\em
$N$-cover} of the graph $\Gamma$, and $N$ the {\em covering
transformation group}. In particular, if $N$ is a cyclic or an
elementary abelian group, then we speak of $\widetilde{\Gamma} $ as
a {\em cyclic cover} or an {\em elementary abelian cover} of
$\Gamma$. Let ${\cal P}:\widetilde{\Gamma} \mapsto \Gamma$ be a
regular covering projection. An automorphism of $\widetilde{\Gamma}$
is said to be {\em fibre-preserving} if it maps a vertex fibre to a
vertex fibre, and all such fibre-preserving automorphisms form a
group, say $F$, called the {\em fibre preserving group}. When $\widetilde{\Gamma}$ is connected, it is easy
to show that $N$ is the kernel of $\Aut(\widetilde{\Gamma})$ acting on the fibres and $F=N_{\Aut(\widetilde{\Gamma})}(N)$. If
$\widetilde{\Gamma}$ is $F$-arc-transitive, we say that
$\widetilde{\Gamma}$ is a {\em symmetric cover} of $\Gamma$.

Two regular covering
projections ${\cal P}:\widetilde{\Gamma} \mapsto \Gamma$ and
${\cal P}^{'}:\widetilde{\Gamma^{'}} \mapsto \Gamma$ of a graph
$\Gamma$ are isomorphic if there exist an automorphism $\alpha
\in \Aut(\Gamma)$ and an isomorphism
$\widetilde{\alpha}:\widetilde{\Gamma} \mapsto
\widetilde{\Gamma^{'}}$ such that $\alpha
{\cal P}={\cal P}^{'}\widetilde{\alpha}$. If $\alpha $ is
identity, then ${\cal P}$ and ${\cal P}^{'}$ are {\em equivalent}, and if $\widetilde{\Gamma}=\widetilde{\Gamma^{'}}$ and
${\cal P}={\cal P}^{'}$ then we call $\widetilde{\alpha}$ a {\em lift} of
$\alpha$ and $\a$ a {\em projection} of $\widetilde{\a}$ along ${\cal P}$.

Let $\Gamma$ be a graph and let $N$ be a finite group. Assign to
each arc $(u,v)$ of $\Gamma$ a voltage $\zeta(u, v) \in N$ such that
$\zeta(u,v)=\zeta(u,v)^{-1}$, where $\z: \G \mapsto N$ is called a
{\em voltage assignment} of $\G$. Let $Cov(\G;\z)$ be the {\em
derived graph} from $\z$, which has vertex set $V(\G) \times N$ and
adjacency relation defined by $(u, a)\thicksim (v, a\z(u,v))$, where
$a\in N$ and $u \thicksim v$ in $\G$. The projection onto the first
coordinate ${\cal P}: Cov(\G; \z) \mapsto \G$ is a regular
$N$-covering projection, where the group $N$ acts semiregularly via
left multiplication on the second coordinate of $(u,a)$, $u\in V(\G)$, $a\in N$. Give a spanning tree $T$ of $\G$, a voltage assignment $\z$ is
said to be {\em T-reduced} if the voltages on the tree arcs are
identity. Gross and Tucker \cite{GT} showed that every regular
covering of a graph $\G$ can be derived from a $T$-reduced voltage
assignment $\z$ with respect to an arbitrary fixed spanning tree of
$\G$. It is clear that if $\z$ is reduced, then the derived graph
$\G\times_{\z} N$ is connected if and only if the voltages on the
cotree arcs generate the voltage group $N$.

A voltage assignment on arcs can be extended to a voltage assignment on walks in a
natural way. Given $\a \in \Aut(\G)$, we define a function $\oa$
from the set of voltages on fundamental closed walks based at a
fixed vertex $v \in V(\G)$ to the voltage group $N$ by $(\z
(C))^{\oa}=\z(C^{\a})$, where $C$ ranges over all fundamental closed
walks at $v$, and $\z(C)$ and $\z(C^{\a})$ are the voltages of $C$
and $C^{\a}$, respectively. Clearly, if $N$ is abelian, then $\oa$ does not depend on the choice of the base vertex, and the fundamental
closed walks at $v$ can be substituted by the fundamental cycles
generated by the cotree arcs of $\G$.  The next proposition is a
special case of Theorem~4.2 in \cite{M}.

\begin{prop}\label{lifts}
Let ${\cal P}: \widetilde{\G}=Cov(\G;\z) \mapsto \G$ be a regular $N$-covering projection. Then an automorphism
$\a$ of $\G$ lifts if and only if $\oa$ extends to an automorphism
of $N$.
\end{prop}

By \cite[Corollary~3.3(a)]{MMP1}, we have the following proposition.

\begin{prop} \label{Iso}
Let ${\cal P}_1: Cov(\G; \z_1) \mapsto \G$ and ${\cal P}_2: Cov(\G;
\z_2) \mapsto \G$ be two regular $N$-covering projections of a graph
$\G$. Then ${\cal P}_1$ and ${\cal P}_2$ are isomorphic if and only
if there is an automorphism $\d \in \Aut(\G)$ and an automorphism
$\eta \in \Aut(N)$ such that $(\z_1(W))^\eta=\z_2(W^\d)$ for all
fundamental closed walks $W$ at some base vertex of $\G$.
\end{prop}

By Theorem~\ref{lem-3.1}~(5), ${\cal CGD}_{p^2}^1$ ($p=5$ or $5\ |\
(p-1)$), ${\cal CGD}_{p^2}^2$ ($5\ |\ (p\pm 1)$), ${\cal CGD}_{p^3}$
($p=5$ or $5\ |\ (p-1)$) and ${\cal CGD}_{p^4}$ are symmetric
elementary abelian covers of  $\Dip_5$. Note that all these graphs have girth $6$.

\begin{theorem} \label{dip5-covers}
Let $p$ be a prime and $\mz_p^n$ an elementary abelian group with
$n\geq 2$. Let $\G$ be a connected symmetric $\mz_p^n$-cover of the dipole
$\Dip_5$. Then $2\leq n\leq 4$ and

\begin{enumerate}

\item [{\rm (1)}] For $n=2$, $\G\cong {\cal CGD}_{p^2}^1$ $(p=5$ or $5\ |\
(p-1))$ or ${\cal CGD}_{p^2}^2$ $(5\ |\ (p\pm 1))$, which are unique
for a given order; $\Aut({\cal CGD}_{5^2}^1)=(R(GD_{5^2})\rtimes
F_{20})\mz_4\cong\mz_5\cdot((F_{20}\times F_{20})\rtimes\mz_2)$ with
$N_A(R(GD_{5^2}))=R(GD_{5^2})\rtimes F_{20}$, $\Aut({\cal
CGD}_{p^2}^1)=R(GD_{p^2})\rtimes\mz_5$ for $5\mid (p-1)$, and
$\Aut({\cal CGD}_{p^2}^2)=R(GD_{p^2})\rtimes D_5$;

\item [{\rm (2)}] For $n=3$, $\G\cong {\cal CGD}_{p^3}$ $(p=5$ or $5\ |\
(p-1))$, which are unique for a given order; $\Aut({\cal
CGD}_{5^3})=R(GD_{5^3})\rtimes S_5$ and  $\Aut({\cal
CGD}_{p^3})=R(GD_{p^3})\rtimes\mz_5$ for  $5\mid (p-1)$;

\item [{\rm (3)}] For $n=4$, $\G\cong {\cal CGD}_{p^4}$ and $\Aut({\cal CGD}_{p^4})=R(GD_{p^4})\rtimes S_5$.
\end{enumerate}
\end{theorem}

\begin{figure}[ht]
\begin{center}
\includegraphics[height=4.5cm,width=5cm]{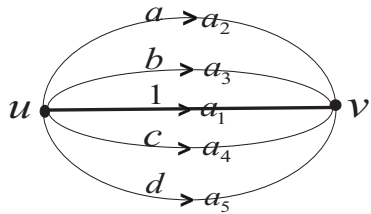}
\end{center}
\vskip -2cm
\caption{\label{picture} $\Dip_{5}$ with a voltage
assignment $\z$.}
\end{figure}

\demo Let $A=\Aut(\G)$. Let $N=\mz_p^n$ and $G=N_A(N)$. Then $N$ has
two orbits and is semiregular on $V(\G)$. By hypothesis, $\G$ is
$G$-arc-transitive. Clearly, $\G_N$ is $\Dip_5$ and its
vertices are denoted by $u$ and $v$ that are connected by five multiple edges
(see Fig.~\ref{picture}).

Let ${\cal P}: \G= \Dip_{5} \times_{\z} \mz_p^n \mapsto \Dip_{5}$ be
the corresponding covering projection, where $\z$ is its voltage
assignment. Since $N\lhd G$, the projection $L$ of $G$ is an
arc-transitive subgroup of $\Aut(\Dip_{5})$. Thus, $L$ lifts to $G$
and $G/N\cong L$. Furthermore, $L$ is the largest subgroup of
$\Aut(\Dip_5)$ which can be lifted along ${\cal P}$. Label the five
arcs of $\Dip_{5}$ starting from $u$  by $a_1$, $a_2$, $a_3$, $a_4$
and $a_5$, respectively. Let $T$ be a spanning tree of $\Dip_5$
corresponding to the arc $a_1$. We may assume that ${\cal P}$ is
$T$-reduced. Write $\z(a_1)=1$, $\z(a_2)=a$, $\z(a_3)=b$,
$\z(a_4)=c$ and $\z(a_5)=d$. Since $\G$ is connected,
$N=\mz_p^n=\langle a, b, c, d\rangle$, forcing that $n\leq 4$.

Let $\alpha =(a_1\ a_2\ a_3\  a_4\  a_5)(a_1^{-1}\  a_2^{-1}\
a_3^{-1}\  a_4^{-1}\ a_5^{-1})$, $\b =( a_1\ a_2\ a_4\ a_3)(
a_1^{-1}\ a_2^{-1}\ a_{4}^{-1}\ a_3^{-1})$, $\gamma= (a_1\
a_1^{-1})(a_2\ a_2^{-1})(a_3\ a_3^{-1})(a_4\ a_4^{-1})(a_5\
a_5^{-1})$,  $\d=(a_1\ a_2\ a_3)(a_1^{-1}\ a_2^{-1}\ a_3^{-1})$, and
$\varepsilon=(a_1\ a_2\ a_5\ a_4)(a_1^{-1}\ a_2^{-1}\ a_5^{-1}\
a_{4}^{-1})$. Then $\a,\b,\g,\d,\varepsilon\in\Aut(\Dip_5)$. There
are four fundamental cycles: $W_1=a_2a_1^{-1}$, $W_2=a_3a_1^{-1}$,
$W_3=a_4a_1^{-1}$ and $W_4=a_5a_1^{-1}$. We list all these cycles
and their voltages in Table~\ref{fig2}, in which $\z(W_i)$ denotes
the voltage on $W_i$.

\begin{table}[ht]
\center
\begin{tabular}{|l|l|l|l|l|l|l|l|}
\hline $W_i$   &$\z(W_i)$ &$W_i^\alpha$
 &$\z(W_i^\alpha)$  &$W_i^\b$ & $\z(W_i^\b)$ & $W_i^{\b^2}$ & $\z(W_i^{\b^2})$

\\\hline
$W_1=a_2a_1^{-1}$   &$a$      &$a_3a_2^{-1}$ &$ba^{-1}$ &$a_4a_2^{-1}$ & $ca^{-1}$ & $a_3a_4^{-1}$ & $bc^{-1}$

\\\hline $W_2=a_3a_1^{-1}$   &$b$ &$a_4a_2^{-1}$ &$ca^{-1}$ &$a_1a_2^{-1}$&$a^{-1}$&$a_2a_4^{-1}$& $ac^{-1}$

\\\hline $W_3=a_4a_1^{-1}$ &$c$ &$ a_5a_2^{-1}$ &$da^{-1}$ &$a_3a_2^{-1}$&$ba^{-1}$&$a_1a_4^{-1}$&$c^{-1}$

\\\hline $W_4=a_5a_1^{-1}$ &$d$ &$
a_1a_2^{-1}$ &$a^{-1}$ &$a_5a_2^{-1}$&$da^{-1}$&$a_5a_4^{-1}$& $dc^{-1}$

\\\hline
\hline

$W_i$   &$\z(W_i)$ & $W_i^\gamma$ & $\z(W_i^\gamma)$
 & $W_i^\d$ & $\z(W_i^\d)$ & $W_i^\varepsilon$
 & $\z(W_i^\varepsilon)$

\\\hline
$W_1=a_2a_1^{-1}$   &$a$        & $a_2^{-1}a_1$ &$a^{-1}$  & $a_3a_2^{-1}$ & $ba^{-1}$& $a_5a_2^{-1}$&$da^{-1}$

\\\hline $W_2=a_3a_1^{-1}$   &$b$  &$a_3^{-1}a_1$ &$b^{-1}$ & $a_1a_2^{-1}$& $a^{-1}$&$a_3a_2^{-1}$ & $ba^{-1}$

\\\hline $W_3=a_4a_1^{-1}$ &$c$ &$ a_4^{-1}a_1$ &$c^{-1}$ &$a_4a_2^{-1}$&$ca^{-1}$& $a_1a_2^{-1}$ & $a^{-1}$

\\\hline $W_4=a_5a_1^{-1}$ &$d$  &$ a_5^{-1}a_1$ &$d^{-1}$ &$a_5a_2^{-1}$&$da^{-1}$& $a_4a_2^{-1}$& $ca^{-1}$

\\\hline

\end{tabular}
\caption{Voltages on fundamental cycles and images under $\a$, $\b$, $\b^2$, $\gamma$, $\d$, $\varepsilon$} \label{fig2}
\end{table}

Clearly, $\Aut(\Dip_{5})=S_5\rtimes \mz_2$, where $S_5$ fixes $u$
and $v$. By the arc-transitivity of $L$, $|L|$ is divisible by $10$.
Let $L^*$ be the subgroup of $L$ fixing $u$ and $v$. Then
$|L:L^*|=2$. Thus, $L^*\leq S_5$ and a Sylow $5$-subgroup of $L^*$
is also a Sylow $5$-subgroup of $\Aut(\Dip_5)$. Since Sylow
$5$-subgroups of $\Aut(\Dip_{5})$ are conjugate, we may assume
$\a\in L$, that is, $\a$ lifts. Noting that $\a^\b=\a^2$, we have that $L^*=\langle \a\rangle, \langle \a,\b^2\rangle,
\langle \a,\b\rangle, A_5$ or $S_5$. In particular, if $\b^2$ cannot
lift then $L^* =\langle \a\rangle\cong\mz_5$; if $\b^2$ lifts but $\b$ and $\d$ cannot then $L^* =\langle \a,\b^2\rangle\cong D_5$;
if $\b$ lifts but $\d$ cannot then $L^*\cong \mz_5\rtimes\mz_4$;
if $\a$, $\b$ and $\d$ lift then $L^*=S_5$.

Consider the mapping $\bar{\a}$ from the set of voltages on the four fundamental cycles of $\Dip_{5}$ to the elementary abelian group $\mz_p^n$, defined by
$\z(W_i)^{\bar{\a}}=\z(W_i^\a)$, $1\leq i\leq 4$. Similarly, we may define $\bar{\b}$, $\bar{\b^2}$, $\bar{\gamma}$, $\bar{\d}$ and $\bar{\varepsilon}$.
By Table \ref{fig2},
$a^{\bar{\gamma}}=a^{-1}$, $b^{\bar{\gamma}}=b^{-1}$,
$c^{\bar{\gamma}}=c^{-1}$, and $d^{\bar{\gamma}}=d^{-1}$. Since $\mz_p^{n}$ is abelian, $\bar{\gamma}$ can be extended to an automorphism of $\mz_p^n$. By Proposition \ref{lifts}, $\g$ lifts along ${\cal P}$. This implies that $L$ lifts if and only if $L^*$ lifts.

Since $\a$ lifts, by Proposition \ref{lifts}, $\bar{\a}$ can be
extended to an automorphism of $\mz_p^n$, denoted by $\a^\ast$.
Again by Table~\ref{fig2},

\parbox{10cm}{
\begin{equation*}
a^{\a^\ast}=ba^{-1}, b^{\a^\ast}=ca^{-1}, c^{\a^\ast}=da^{-1},
d^{\a^\ast}=a^{-1}.
\end{equation*}}\hfill
\parbox{1cm}{
\begin{equation}
\label{eq0} \end{equation}}

Suppose $b=1$. By Eq~(\ref{eq0}), $ca^{-1}=b^{\a^\ast}=1$, that is, $c=a$. Thus, $d=b=1$ because $c^{\a^*}=a^{\a^*}$, and hence $c^{-1}=a^{-1}=d^{\a^*}=1$. It follows that
$\mz_p^{n}= \langle a, b, c, d\rangle=1$, a contradiction. Thus, $b\neq 1$. Similarly,  $a \neq 1$, $c \neq 1$ and $d \neq 1$.
Since $n \leq 4$, we have $N=\mz_p^2$, $\mz_p^3$ or $\mz_p^4$.

Suppose $d\in \langle a \rangle$. Then $d=a^k$, $k \in \mz_p^{\ast}$. By Eq~(\ref{eq0}), $a^{-1}=b^{k}a^{-k}$, implying
$b \in \langle a \rangle$. Similarly, $c \in \langle a \rangle$. It follows that $\mz_p^{2}= \langle a, b, c, d\rangle=\langle a
\rangle$, a contradiction. Thus,  $d \notin \langle a \rangle$.

\medskip
\f \f {\bf Case 1:}. $N=\mz_p^2$.

Since $d \notin \langle a \rangle$, we have $\mz_p^2=\langle a, d\rangle$. Let $b=a^id^j$ and $c=a^{\ell}d^m$ for some $i,j,\ell, m \in \mz_p $. By Eq~(\ref{eq0}),
$a^{i(i-1)-j}d^{ij}=(a^{i-1}d^j)^ia^{-j}=(ba^{-1})^ia^{-j}=
(a^id^j)^{\a^\ast}=b^{\a^\ast}=ca^{-1}=a^{\ell-1}d^{m}$
and $a^{(i-1)\ell-m}d^{j\ell}=(a^{i-1}d^j)^\ell a^{-m}=(ba^{-1})^\ell a^{-m}=(a^\ell
d^m)^{\a^\ast}=c^{\a^\ast}=a^{-1}d$. Then the following equations hold in the field $\mz_p$.

\parbox{6cm}{
\begin{eqnarray*}
&& i(i-1)-j=\ell-1,\\
&& ij=m,\\
&& (i-1)\ell-m=-1,\\
&& j\ell=1.\\
\end{eqnarray*}}\hfill
\parbox{1.0cm}{
\begin{eqnarray}
\label{eq1} \\ \label{eq2}\\ \label{eq3}\\
\label{eq4}\end{eqnarray}}

By Eqs~(\ref{eq4}) and (\ref{eq2}), $\ell\neq  0$, $j=\ell^{-1}$
and $i=m\ell$. By Eq~(\ref{eq3}), $\ell(m\ell-1)-m=-1$, that is,
$(\ell-1)(m\ell+m-1)=0$. Thus, $\ell-1=0$ or $m\ell+m-1=0$ in $\mz_p$.

Assume $\ell=1$. Then $j=\ell^{-1}=1$ and $i=m\ell=m$. By
Eq~(\ref{eq1}), $i(i-1)=1$, which implies that $p\not=2$ and $(2i-1)^2=5$. By elementary
number theory, $p=5k\pm 1$ or $p=5$.

Let $p=5k\pm 1$. Since $(2i-1)^2=5$, we have $i=2^{-1}(1+\lambda)$,
where $\ld^2=5$. It follows that $b=a^{2^{-1}(1+\lambda)}d$ and
$c=ad^{2^{-1}(1+ \lambda)}$. Note that
$2^{-1}(1+\lambda)-1=-2^{-1}(1-\lambda)$. By Table~\ref{fig2},
$\bar{\a}$ can be extended to the automorphism of $\mz_p^2$ induced
by $a\mapsto a^{-2^{-1}(1-\lambda)}d$ and $d\mapsto a^{-1}$, and
$\bar{\b^2}$ can be extended to the automorphism induced by
$a\mapsto a^{-2^{-1}(1-\ld)}d^{2^{-1}(1-\ld)}$ and $d\mapsto
a^{-1}d^{2^{-1}(1-\ld)}$. Furthermore, $\bar{\b}$ and $\bar{\d}$
cannot be extended to automorphisms of $\mz_p^2$. Thus, $L^*\cong
D_5$ and $|G|=|N||L|=20p^2$.

We claim that $\G$ is unique for any prime $p$ such that $p=5k\pm
1$. This is sufficient to show that  $\G$ is independent of the
choice of $\ld$. Note that the equation $x^2=5$ has exactly two roots
in $\mz_p$, that is, $\pm \ld$. Since $(2i-1)^2=5$, we have
$i=2^{-1}(1\pm \lambda)$. It follows that $b=a^{2^{-1}(1\pm
\lambda)}d$ and $c=ad^{2^{-1}(1\pm \lambda)}$.  Clearly, the voltage
assignment $\z$ is determined by $\z(a_1)$, $\z(a_2)$, $\z(a_3)$,
$\z(a_4)$ and $\z(a_5)$, and for convenience, write $\z=(\z(a_1),
\z(a_2), \z(a_3), \z(a_4), \z(a_5))$. It follows that  $\z=\z_1=(1, a,
a^{2^{-1}(1+\ld)}d, ad^{2^{-1}(1+\ld)}, d)$ or $\z=\z_2=(1, a,
a^{2^{-1}(1-\ld)}d, ad^{2^{-1}(1-\ld)}, d)$. By Table~\ref{fig2},
$\z_1(W_1)=a$, $\z_1(W_2)=a^{2^{-1}(1+\ld)}d$,
$\z_1(W_3)=ad^{2^{-1}(1+\ld)}$, $\z_1(W_4)=d$ and
$\z_2(W_1^\b)=d^{2^{-1}(1-\ld)}$, $\z_2(W_2^\b)=a^{-1}$,
$\z_2(W_3^\b)=a^{2^{-1}(1-\ld)-1}d$, $\z_2(W_4^\b)=da^{-1}$. Let
$\eta$ be the automorphism of $\mz_p^2$ induced by $a \mapsto
d^{2^{-1}(1-\lambda)}$ and $d \mapsto da^{-1}$. Then
$(\z_1(W_i))^\eta=\z_2(W_i^\b)$ for each $1\leq i\leq 4$, and by
Proposition \ref{Iso}, $\Dip_5 \times_{\z_{1}} \mz_p^{2}\cong\Dip_5
\times_{\z_{2}} \mz_p^{2}$, as claimed.

Let $p=5$. Then $(2i-1)^2=0$ implies $i=3$ in $\mz_5$. Thus,
$m=i=3$. It follows that $b=a^{3}d$ and $c=ad^{3}$. By
Table~\ref{fig2}, $\bar{\a}$ and $\bar{\b}$ can be extended to the
automorphisms of $\mz_p^2$ induced by $a\mapsto a^2d$, $d\mapsto
a^{-1}$, and $a\mapsto d^2$, $d\mapsto da^{-1}$, respectively. But,
$\bar{\d}$ can not be extended to an automorphism of $\mz_p^2$.
Thus, $L^*\cong \mz_5\rtimes\mz_4$ and $|G|=40\times 5^2$.

Assume $m\ell+m-1=0$. Then $\ell+1\not=0$ and $m=(\ell+1)^{-1}$.
Recall that $j=\ell^{-1}$ and  $i=m\ell$. Then
$i=\ell(\ell+1)^{-1}$. By Eq~(\ref{eq1}),
$\ell(\ell+1)^{-1}(\ell(\ell+1)^{-1}-1)-\ell^{-1}=\ell-1$, that is,
$\ell^4+\ell^3+\ell^2+\ell+1=0$. Thus,
$\ell^5-1=(\ell-1)(\ell^4+\ell^3+\ell^2+\ell+1)=0$. For $p=5$, we
have $\ell=1$ because $\ell^5=1$ in $\mz_5^*$, which has been
discussed in the previous paragraph. Let $p>5$. Then $\ell$ is an
element of order $5$ in $\mz_p^*$ and hence $5\ |\ (p-1)$. In this
case, $(\ell+1)^{-1}=-\ell^3-\ell$ and
$(\ell^2+1)^{-1}=-\ell^2-\ell$. Note that
$b=a^{\ell(\ell+1)^{-1}}d^{\ell^{-1}}$ and
$c=a^{\ell}d^{(\ell+1)^{-1}}$. By Table~\ref{fig2}, $\bar{\a}$ can
be extended to the automorphism of $\mz_p^2$ induced by $a\mapsto
a^{\ell(\ell+1)^{-1}-1}d^{\ell^{-1}}$ and $d\mapsto a^{-1}$, but
$\bar\b^2$ cannot. Thus, $L^*\cong \mz_5$ and $|G|=10p^2$.

We now claim that $\G$ is unique for any prime $p$ with $5\ |\
(p-1)$. Clearly, $\mz_p^*$ has exactly four elements of order $5$,
that is, $\ell^i$ for $1\leq i\leq 4$. By the arbitrariness of $\ell$,
it suffices to show that $\Dip_5 \times_{\z_{1}} \mz_p^{2}\cong
\Dip_5 \times_{\z_{2}} \mz_p^{2}$, where $\z_1$ and $\z_2$ are
voltage assignments corresponding to $\ell$ and $\ell^2$
respectively, that is, $\z_1=(1, a,
a^{\ell(\ell+1)^{-1}}d^{\ell^{-1}}, a^{\ell}d^{(\ell+1)^{-1}}, d)$
and $\z_2=(1, a, a^{\ell^{2}(\ell^{2}+1)^{-1}}d^{\ell^{-2}},
a^{\ell^{2}}d^{(\ell^{2}+1)^{-1}}, d)$. Recall that
$(\ell+1)^{-1}=-\ell^3-\ell$ and
$(\ell^2+1)^{-1}=-\ell^2-\ell$.
By Table~\ref{fig2},
$\z_1(W_1)=a$, $\z_1(W_2)=a^{\ell^{3}+\ell+1}d^{\ell^{4}}$,
$\z_1(W_3)=a^{\ell}d^{\ell^{4}+\ell^{2}+1}$, $\z_1(W_4)=d$, and
$\z_2(W_1^{\varepsilon})=a^{-1}d$,
$\z_2(W_2^{\varepsilon})=a^{\ell^{2}+\ell}d^{\ell^{3}}$,
$\z_2(W_3^{\varepsilon})=a^{-1}$, $\z_2(W_4^{\varepsilon})=
a^{\ell^{2}-1}d^{\ell^{4}+\ell^{3}+1}$. Let $\eta$ be the
automorphism of $\mz_p^2$ induced by $a \mapsto a^{-1}d$ and $d
\mapsto a^{\ell^{2}-1}d^{\ell^{4}+\ell^{3}+1}$. It is easy to check
that $(\z_1(W_i))^\eta=\z_2(W_i^\varepsilon)$ for each $1\leq i\leq
4$, and by Proposition~\ref{Iso}, $\Dip_5 \times_{\z_{1}}
\mz_p^{2}\cong \Dip_5 \times_{\z_{2}} \mz_p^{2}$, as claimed.

We have proved that for $p=5$, the graph $\G$ is unique with
$|N_A(N)|=40\times 5^2$; for $5\ |\ (p-1)$, the graph $\G$ is unique
with $|N_A(N)|=10p^2$; for $5\ |\ (p\pm 1)$, the graph $\G$ is
unique with $|N_A(N)|=20p^2$. By Theorem~\ref{lem-3.1}~(1) and (2), the uniqueness of $\G$ implies that $\G\cong {\cal CGD}_{p^2}^1$ ($p=5$ or $5\ |\ (p-1))$ or ${\cal
CGD}_{p^2}^2$ ($5\ |\ (p\pm 1))$.

Let $\G\cong {\cal CGD}_{p^2}^1$ ($p=5$ or $5\ |\ (p-1))$.  By
MAGMA~\cite{BCP}, $|\Aut({\cal CGD}_{5^2}^1)|=2^5\cdot5^3$ and
$\Aut({\cal CGD}_{5^2}^1)$ has a normal subgroup of order $5$. By
Theorem~\ref{lem-3.1} and Proposition~\ref{stabilizer-val5},
$N_A(R(GD_{5^2}))=R(GD_{5^2})\rtimes F_{20}$ and
$A=(R(GD_{5^2})\rtimes F_{20})\mz_4\cong\mz_5\cdot((F_{20}\times
F_{20})\rtimes\mz_2)$, where $(F_{20}\times F_{20})\rtimes\mz_2$ is
a maximal solvable subgroup of $\Aut(K_{5,5})\cong (S_5\times
S_5)\rtimes\mz_2$. For $5\ |\ (p-1)$, we claim that $|\Aut({\cal
CGD}_{p^2}^1)|=10p^2$. Again by MAGMA~\cite{BCP}, $|\Aut({\cal
CGD}_{p^2}^1)|=10p^2$ for $p=11, 31, 41, 61, 71$. Assume $p\geq 91$.
Since $|N_A(N)|=10p^2$, it suffices to show $N\lhd A$. Suppose
$N\ntriangleleft A$. Since $A=R(GD_{p^2})A_w$ for any $w\in V(\G)$,
$N$ is a Sylow $p$-subgroup of $A$. By Sylow Theorem, the number of
Sylow $p$-subgroups of $A$ is $kp+1$ and $kp+1=|A:N_A(N)|$. It follows that $kp+1=|A_w|/5$ because
$N_A(N)=R(GD_{p^2})\rtimes \mz_5$ by Theorem~\ref{lem-3.1}. Since ${\cal CGD}_{p^2}^1$ has girth $6$, $\G$ is at
most $3$-arc-transitive and by Proposition~\ref{stabilizer-val5},
$|A_w|\ |\ 2^6\cdot 3^2\cdot 5$. Thus, $kp+1$ is a divisor of
$2^6\cdot 3^2$. Since $p\geq 91$ and $5\ |\ (p-1)$, we have
$(k,p)=(1,191)$ and $kp=2^6\cdot 3-1$. However, this is impossible
by Proposition~\ref{stabilizer-val5} because $|A_w|=5(kp+1)=2^6\cdot
3\cdot 5$. Thus, $|\Aut({\cal CGD}_{p^2}^1)|=10p^2$ and
$A=R(GD_{p^2})\rtimes\mz_5$.

Let $\G\cong {\cal CGD}_{p^2}^2$ ($5\ |\ (p\pm 1))$.  We claim
$|\Aut({\cal CGD}_{p^2}^2)|=20p^2$. By MAGMA~\cite{BCP}, this is
true for $p=11, 19,29,31,41,59,61,71$. Assume $p\geq 79$ and it
suffices to show that $N\lhd A$ because $N_A(N)=R(GD_{p^2})\rtimes
D_5$. Suppose $N\ntriangleleft A$. As the previous paragraph, the
number  $kp+1$ of Sylow $p$-subgroups of $A$ is a divisor of
$|A_w|/10=2^5\cdot 3^2$. However, no such a prime $p$ exists for
$p\geq 79$. Thus, $|\Aut({\cal CGD}_{p^2}^2)|=20p^2$ and
$A=R(GD_{p^2})\rtimes D_5$.

\medskip
\f \f {\bf Case 2:} $N=\mz_p^3$.

Recall that $d\not\in \langle a \rangle$.
Suppose $b\in \langle a, d \rangle$. Then $b=a^{i}d^{j}$, $i, j \in \mz_p$. By Eq~(\ref{eq0}),
$ca^{-1}=b^{\a^*}=(b^ia^{-i})a^{-j}$, implying $c\in\langle a, d \rangle$. It
follows that $\mz_p^{3}=\langle a, b, c, d\rangle=\langle a,
d\rangle$, a contraction. Thus, $b\notin \langle a, d \rangle$ and
$\mz_p^{3}=\langle a, d, b \rangle$. Let $c=a^ib^jd^k$ for some $i,
j, k \in \mz_p$. Again by Eq~(\ref{eq0}),
$da^{-1}=c^{\a^{\ast}}=(a^ib^jd^k)^{\a^{\ast}}=
(ba^{-1})^i(a^{i-1}b^jd^k)^j(a^{-1})^k=
a^{ij-i-j-k}b^{i+j^2}d^{kj}$, and the following equations hold in
the field $\mz_p$.

\parbox{6cm}{
\begin{eqnarray*}
&& ij-i-j-k=-1,\\
&& i+j^2=0,\\
&& kj=1.\\
\end{eqnarray*}}\hfill
\parbox{1.0cm}{
\begin{eqnarray}
\label{eq5} \\ \label{eq6}\\ \label{eq7}
\end{eqnarray}}

By Eqs~(\ref{eq6}) and (\ref{eq7}),  $i=-j^2$ and $k=j^{-1}$. By
Eq~(\ref{eq5}),  $-j^3+j^2-j-j^{-1}+1=0$, that is,
$j^4-j^3+j^2-j+1=0$. Let $\ell=-j$. Then
$\ell^4+\ell^3+\ell^2+\ell+1=0$ and
$\ell^5-1=(\ell-1)(\ell^4+\ell^3+\ell^2+\ell+1)=0$. It follows that
either $p=5$ with $\ell=1$, or $5\ |\ (p-1)$ with $\ell$ an element
of order $5$ in $\mz_p^*$.

Let $p=5$. Then $c=a^{4}b^{4}d^{4}$. By Table~\ref{fig2}, $\bar{\a}$
can be extended to the automorphism of $\mz_p^3$ induced by $a\mapsto
a^4b$, $b\mapsto a^3b^4d^4$, and $d\mapsto a^4$. Thus, $\G$ is
symmetric.

Let $5\ |\ (p-1)$ with $\ell$ an element of order $5$ in $\mz_p^*$.
Then $c=a^{-\ell^2}b^{-\ell}d^{-\ell^4}$. By Table~\ref{fig2},
$\bar{\a}$ can be extended to the automorphism of $\mz_p^3$ induced
by $a\mapsto a^{-1}b$, $b\mapsto a^{-\ell^2-1}b^{-\ell}d^{-\ell^4}$,
and $d\mapsto a^{-1}$. Thus, $\G$ is symmetric. Further,
$\bar{\b^2}$ cannot be extended to an automorphism of $\mz_p^3$. Thus,
$L^*\cong \mz_5$ and $|G|=10p^3$. We now claim that $\G$ is unique.
Since $\mz_p^*$ has exactly four elements of order $5$, that is,
$\ell^i$ for each $1\leq i\leq 4$, it suffices to show that
$\Dip_5 \times_{\z_{1}} \mz_p^{2}\cong \Dip_5 \times_{\z_{2}}
\mz_p^{2}$, where $\z_1$ and $\z_2$ are voltage assignments
corresponding to $\ell$ and $\ell^2$ respectively, that is, $\z_1=(1, a,
b, a^{-\ell^2}b^{-\ell}d^{-\ell^4}, d)$ and $\z_2=(1, a, b,
a^{-\ell^4}b^{-\ell^2}d^{-\ell^3}, d)$. By Table~\ref{fig2},
$\z_1(W_1)=a$, $\z_1(W_2)=b$,
$\z_1(W_3)=a^{-\ell^2}b^{-\ell}d^{-\ell^4}$, $\z_1(W_4)=d$, and
$\z_2(W_1^{\varepsilon})=a^{-1}d$,
$\z_2(W_2^{\varepsilon})=a^{-1}b$, $\z_2(W_3^{\varepsilon})=a^{-1}$,
$\z_2(W_4^{\varepsilon})=a^{-1-\ell^4}b^{-\ell^2}d^{-\ell^3}$. Let
$\eta$ be the automorphism of $\mz_p^2$ induced by $a \mapsto
a^{-1}d$, $b \mapsto a^{-1}b, d \mapsto
a^{-1-\ell^4}b^{-\ell^2}d^{-\ell^{-2}}$. Then
$(\z_1(W_i))^\eta=\z_2(W_i^\varepsilon)$ for each $1\leq i\leq 4$,
and by Proposition~\ref{Iso}, $\Dip_5 \times_{\z_{1}} \mz_p^{2}\cong
\Dip_5 \times_{\z_{2}} \mz_p^{2}$, as claimed.

We have proved that $p=5$ or  $5\ |\ (p-1)$. In each case, the graph
$\G$ is unique and for the latter, $|N_A(N)|=10 p^3$. Then
Theorem~\ref{lem-3.1}~(3) implies that $\G\cong {\cal CGD}_{p^3}$.
By MAGMA~\cite{BCP}, $|\Aut({\cal CGD}_{5^3})|=2^4\cdot 3\cdot 5^4$
and $|\Aut({\cal CGD}_{11^3})|=2\cdot 5\cdot 11^3$, and again by
Theorem~\ref{lem-3.1}~(3), $\Aut({\cal
CGD}_{5^3})=R(GD_{5^3})\rtimes S_5$ and  $\Aut({\cal
CGD}_{11^3})=R(GD_{11^3})\rtimes \mz_5$. Assume $p\geq 31$. We claim
that $|\Aut({\cal CGD}_{p^3})|=10p^3$. It suffices to show $N\lhd A$
because $|N_A(N)|=10p^3$. Suppose to the contrary that
$N\ntriangleleft A$. Since $A=R(GD_{p^3})A_w$ for any $w\in V(\G)$,
$N$ is a Sylow $p$-subgroup of $A$. By Sylow Theorem, the number of
Sylow $p$-subgroups in $A$ is $kp+1$ and $kp+1=|A:N_A(N)|$. Since
$N_A(N)=R(GD_{p^3})\rtimes \mz_5$, we have $kp+1=|A_w|/5$, and since
${\cal CGD}_{p^3}$ has girth $6$, it is at most $3$-arc-transitive.
By Proposition~\ref{stabilizer-val5}, $|A_w|\ |\ 2^6\cdot 3^2\cdot
5$ and hence $kp+1$ is a divisor of $2^6\cdot 3^2$. Since $p\geq 31$
and $5\ |\ (p-1)$, we have $(k,p)=(1,31)$ and $kp=2^5-1$,
$(k,p)=(1,191)$ and $kp=2^6\cdot 3-1$, $(k,p)=(1,71)$ and
$kp=2^3\cdot 3^2-1$, or $(k,p)=(7,41)$ and $kp=2^5\cdot 3^2-1$. It
follows that $|A_w|=2^5\cdot 5$, $2^6\cdot 3\cdot 5$, $2^3\cdot
3^2\cdot 5$ or $2^5\cdot 3^2\cdot 5$. Since $\G$ is not $4$-arc-transitive, by
Proposition~\ref{stabilizer-val5} only the last case can happen and
in this case,  $p=41$ and $\G\cong {\cal CGD}_{41^3}$ is
$3$-arc-transitive. However, it is easy to check by
MAGMA~\cite{BCP} that there is no $6$-cycle in ${\cal CGD}_{41^3}$
passing through the $3$-arc $(1,h,ah\cdot h,h\cdot ah\cdot
h)=(1,h,a,a^{-1}h)$ (one may take $\ell=10$ because $10^5=1$ in
$\mz_{41}$), which is impossible because of the
$3$-arc-transitivity of $\G$. Thus, $|\Aut({\cal CGD}_{p^3})|=10p^3$
for $5\ |\ (p-1)$ and $\Aut({\cal CGD}_{p^3})=R(GD_{p^3})\rtimes
\mz_5$.

\medskip

\f \f {\bf Case 3:} $N=\mz_p^4$.

Let $N=\langle a\rangle \times \langle b\rangle \times \langle
c\rangle \times \langle d\rangle$. By the connectedness of $\G$, we
may let $\z(a_1)=1$, $\z(a_2)=a$, $\z(a_3)=b$, $\z(a_4)=c$, and
$\z(a_5)=d$. By Table~\ref{fig2}, it is easy to see that $\bar{\a}$,
$\bar{\b}$ and $\bar{\d}$ can be extended to automorphisms of
$\mz_p^4$. Thus, $L^*=S_5$ and $|G|=240p^4$. By
Theorem~\ref{lem-3.1}, $\G\cong {\cal CGD}_{p^4}$.

For $p=2,3,5$ or $7$, by MAGMA~\cite{BCP} we have $|\Aut({\cal
CGD}_{p^4})|=2^4\cdot 3\cdot 5\cdot p^4$. This is also true for any
prime $p\geq 11$. To prove it, we only need to show that $N\lhd A$
because $|G|=|N_A(N)|=2^4\cdot 3\cdot 5\cdot p^4$. Suppose to the
contrary that $N\ntriangleleft A$. Note that $N$ is a Sylow
$p$-subgroup because $A=R(GD_{p^2})A_w$ for any $w\in V(\G)$. By
Sylow Theorem, the number of Sylow $p$-subgroups in $A$ is $kp+1$
and $kp+1=|A:N_A(N)|$. Since $N_A(N)=R(GD_{p^2})\rtimes S_5$, we
have $kp+1=|A_w|/120$, and since ${\cal CGD}_{p^3}$ has girth $6$,
it is at most $3$-arc-transitive. By
Proposition~\ref{stabilizer-val5}, $|A_w|\ |\ 2^6\cdot 3^2\cdot 5$
and hence $kp+1$ is a divisor of $2^3\cdot 3$. Since $p\geq 11$, we
have $(k,p)=(1,11)$ and $kp=2^2\cdot 3-1$, or $(k,p)=(1,23)$ and
$kp=2^3\cdot 3-1$. It follows that $|A_w|=2^5\cdot 3^2\cdot 5$ or
$2^6\cdot 3^2\cdot 5$, and by Proposition~\ref{stabilizer-val5},
$\G$ is $3$-arc-transitive. Note that $\G\cong {\cal CGD}_{11^4}$ or
${\cal CGD}_{23^4}$. It is easy to check that  there is no $6$-cycle
in ${\cal CGD}_{11^4}$ and  ${\cal CGD}_{23^4}$ passing through the
$3$-arc $(1,h,ah\cdot h,h\cdot ah\cdot h)=(1,h,a,a^{-1}h)$, which
contradicts the $3$-arc-transitivity of $\G$. Thus, $|\Aut({\cal
CGD}_{p^4})|=240p^4$ for each prime $p$, and by
Theorem~\ref{lem-3.1}~(4), $\Aut({\cal
CGD}_{p^4})=R(GD_{p^4})\rtimes S_5$.\hfill\qed

\section{Basic graphs\label{s4}}

In this section we investigate pentavalent symmetric  graphs of order $2p^n$ and basic ones of such graphs are determined. We first prove the following lemma.

\begin{lemma}\label{lem-3.2}
Let $p$ be a prime and let $\G$ be a connected pentavalent symmetric
graph of order $2p^n$ with $n\geq 2$. Let $G$ be an arc-transitive
subgroup of $ \Aut(\G)$. Then every minimal normal subgroup of $G$
is an elementary abelian $p$-group.
\end{lemma}

\demo  Let $u\in V(\G)$. By Proposition~\ref{stabilizer-val5},
$|G_u| \mid 2^9\cdot3^2\cdot 5$, and hence $|G|\ |\
2^{10}\cdot3^2\cdot5\cdot p^n$. Let $N$ be a minimal normal subgroup
of $G$. Then $N\cong T^{m}$ for a finite simple group $T$.

Suppose that $N$ is nonabelian. If $N$ has more than two orbits, by
Proposition~\ref{QG}, $|N|\mid 2p^n$ and hence $N$ is solvable, a
contradiction. Thus, $N$ has one or two orbits on $V(\G)$, which implies that $p^n\mid |N|$. Since
$|G_u|\mid 2^9\cdot 3^2\cdot 5$ and $|V(\G)|=2p^n$, we have
$|G|=2^i\cdot3^j\cdot5\cdot p^n$ with $1\leq i\leq 10$ and $0\leq
j\leq 2$. Thus, $\pi(N)=4$ or $3$, where $\pi(N)$ is the number of
distinct prime factors of $|N|$.

If $\pi(N)=4$ then $p\geq 7$ and since $5^2\nmid |G|$, $N=T$ is a simple
$\{2,3,5,p\}$-group such that $2^{11}\nmid |N|$, $3^3\nmid |N|$. Suppose that
$p^2\mid |N|$. By Proposition~\ref{non-abelian-simple-groups}, we have $5^2
\mid |N|$ and hence $5^2\mid |G|$, a contradiction. Thus, $p^2\nmid
|N|$, and since $p^n\mid |N|$, we have $n=1$, contrary to
the hypothesis. It follows that $\pi(N)=3$. Since $\G$ has the prime valency $5$, $G_{u}$ is primitive on $\Gamma_1(u)$,
the neighborhood of $u$ in $\G$, and since $N_{u}\lhd G_{u}$, either
$N_{u}=1$ or $5\ |\ |N_{u}|$. If $N_u=1$ then $N$ is semiregular on $V(\G)$. Thus, $|N|\ |\ 2p^n$ and hence $N$ is solvable, a contradiction. Thus, $5\ |\ |N_u|$, and by Proposition~\ref{non-abelian-simple-groups}, $N$ is a $\{2,3,5\}$-group. Furthermore, $p=5,2$ or $3$, and $T\cong A_5, A_6$ or $\PSU(4,2)$.

For $p=5$, we have $N=T$ or $T^2$ as $3^3 \nmid |G|$. It follows that $5^3 \nmid |N|$, and since $N$ has at most two orbits and $5\ |\ |N_u|$, we have $n=1$, contrary to the hypothesis.

For $p=2$ or $3$, we have $N=T$ as $5^2\nmid |G|$. If $p=2$ then $N=A_5$ or $A_6$ as $3^3 \nmid |G|$, which implies that $|\G|=8$ or $16$ because $2^n\mid |N|$. By McKay~\cite{Mckay}, there is no pentavalent symmetric  graph of order $8$
and there is a unique pentavalent symmetric graph of order $16$ that is the graph $FQ_4$, but this is impossible because
$\Aut(FQ_4)\cong \mz_2^4\rtimes S_5$ has no normal subgroup $A_5$ or
$A_6$.
Thus, $p=3$.  Since $|\G|=2\cdot 3^n\geq 18$ and since there is no pentavalent symmetric graph of order
$18$ by McKay~\cite{Mckay}, we have $N=$PSU$(4,2)$, and since $|$PSU$(4,2)|=2^6\cdot 3^4\cdot5$, we have $|\G|=2\cdot 3^3$ or $2\cdot 3^4$. If $N$ is transitive on $V(\G)$ then $|N_u|=2^5\cdot 3\cdot 5$ or $2^5\cdot  5$, which is impossible by Proposition~\ref{stabilizer-val5}. Therefore, $N$ has two orbits and $|N_u|=2^6\cdot 3\cdot 5$ or $2^6\cdot  5$.
Clearly, $\G$ is bipartite and $G$ has a $2$-element, say $g$, interchanging the two bipartition sets of $\G$. It follows that $H=N\lg g\rg$ is arc-transitive and $|H_u|=2^{i+6}\cdot 3\cdot 5$, or $2^{i+6}\cdot 5$, $i\geq 0$. Again by Proposition~\ref{stabilizer-val5}, $H_u\cong \ASL(2,4)$ or $\ASIL(2,4)$, and since $N_u\lhd H_u$, we have $N_u\cong \ASL(2,4)$. However, $\PSU(4,2)$ has no subgroup isomorphic to $\ASL(2,4)$ by MAGMA~\cite{BCP}, a contradiction.

Now $N$ is abelian and hence elementary abelian. Since $|V(\G)|=2p^n$, $N$ is a $2$-group or a
$p$-group. If $p\not=2$ and $N$ is a $2$-group then the quotient
graph $\G_N$ has odd order and valency $5$, a contradiction. Thus,
$N$ is a $p$-group. This completes the proof. \qed

The {\em hypercube $Q_n$}
is the Cayley graph $\Cay(\mz_2^n,\{a_1,a_2,\cdots,a_n\})$, where
$\mz_2^n=\langle a_1\rangle\times\langle
a_2\rangle\times\cdots\times\langle a_n\rangle$, and it
is well-known that $\Aut(Q_n)=R(\mz_2^n)\rtimes S_n$. The {\em folded
hypercube $FQ_n$} $(n\geq 2)$ is the Cayley graph
$FQ_{n}=\Cay(\mz_2^n, S)$, where $S=\{a_1, a_2, \cdots, a_n,$
$a_1a_2\cdots a_n\}$. It is easy to see that $\Aut(\mz_2^n,S)$ is
transitive on $S$ and any permutation on $\{a_1, a_2, \cdots,a_n\}$
induces an automorphism of $\mz_2^n$, which fixes $S$ setwise. This
implies that $S_{n+1}\leq\Aut(\mz_2^n,S)$. For $n\geq 4$, $FQ_n$ has a unique $4$-cycle passing through $1$ and any two elements in $S$, and hence $\Aut(FQ_n)=R(\mz_2^n)\rtimes
S_{n+1}$.

\begin{theorem}\label{th-3.1}
Let $p$ be a prime and let $\G$ be a connected pentavalent symmetric
graph of order $2p^n$ with $n\geq 1$. Then $\G$ is a normal cover of one of
the following graphs: $K_6$,  $FQ_4$, ${\cal CD}_{p}$ $(p=5$ or $5\
|\ (p-1))$, ${\cal CGD}_{p^2}^1$ $(p=5$ or $5\ |\ (p-1))$, ${\cal
CGD}_{p^2}^2$ $(5\ |\ (p\pm 1))$, ${\cal CGD}_{p^3}$ $(p=5$ or $5\
|\ (p-1))$, or ${\cal CGD}_{p^4}$.
\end{theorem}

\demo Let $A=\Aut(\G)$ and let
$M$ be a maximal normal subgroup of $A$ which has more than two
orbits on $V(\G)$. By Proposition~\ref{QG}, $\G$ is a normal cover
of the quotient graph $\G_M$, which has valency $5$ and is
$A/M$-arc-transitive. Clearly, $|V(\G_M)|=2p^m$ for a positive integer
$m$. We aim to show that $\G_M\cong K_6$,  $FQ_4$, ${\cal CD}_{p}$
$(p=5$ or $5\ |\ (p-1))$, ${\cal CGD}_{p^2}^1$ $(p=5$ or $5\ |\
(p-1))$, ${\cal CGD}_{p^2}^2$ $(5\ |\ (p\pm 1))$, ${\cal CGD}_{p^3}$
$(p=5$ or $5\ |\ (p-1))$, or ${\cal CGD}_{p^4}$.

Let $N$ be a minimal normal subgroup of $A/M$. If $m=1$ then, by
Proposition~\ref{2p}, $\G_M\cong K_6$, or ${\cal CD}_{p}$ $(p=5$ or
$5\ |\ (p-1))$. In what follows we assume that $m\geq 2$. By
Lemma~\ref{lem-3.2}, $N$ is an elementary abelian $p$-group, and by the
maximality of $M$, $N$ has one or two orbits on $V(\G_M)$.

First assume that $N$ has one orbit, that is, $N$ is transitive on
$V(\G_M)$. Then $N$ acts regularly on $V(\G_M)$ and $|N|=|V(\G_M)|$.
It follows that $N$ is an elementary abelian $2$-group and $\G_M=\Cay(N,S)$.
Since $\G_M$ has valency $5$, the connectedness of $\G_M$ implies
that $|N|\leq 32$. Clearly, if $|N|=32$ then $\G_M\cong Q_5={\cal
CGD}_{2^4}$. By McKay~\cite{Mckay}, there is no pentavalent
symmetric  graph of order $8$ and there is a unique pentavalent
symmetric  graph of order $16$. Thus, if $|N|=16$ then $\G_M\cong
FQ_4$.

Now assume that $N$ has two orbits on $V(\G_M)$. Then $\G_M$ is a
bipartite graph with the two orbits of $N$ as its bipartite sets. Let
$u\in V(\G_M)$. If $N_u\not=1$ then $\G_M\cong K_{5,5}$, contrary to the assumption $m\geq 2$. Thus, $N$ is regular
on each bipartite set of $\G_M$ and hence $\G_M$ is a symmetric $N$-cover
of the dipole $\Dip_5$. By Theorem~\ref{dip5-covers}, $\G_M\cong {\cal
CGD}_{p^2}^1$ $(p=5$ or $5\ |\ (p-1))$, ${\cal CGD}_{p^2}^2$ $(5\ |\
(p\pm 1))$, ${\cal CGD}_{p^3}$ $(p=5$ or $5\ |\ (p-1))$, or ${\cal
CGD}_{p^4}$. \hfill\qed

\begin{theorem}\label{Basic}
Let $p$ be a prime and $n$ a positive integer. The basic graphs of connected pentavalent symmetric
graph of order $2p^n$ are $K_6$, $FQ_4$, ${\cal CGD}_{5^3}$, ${\cal
CD}_{p}$ $(p=5$ or $5\ |\ (p-1))$, ${\cal CGD}_{p^2}^2$ $(5\ |\
(p\pm 1))$ or ${\cal CGD}_{p^4}$ $(p=3$ or $p\geq 7)$.
\end{theorem}

\demo Let $\G=K_6$, $FQ_4$, ${\cal CGD}_{5^3}$, ${\cal
CD}_{p}$ $(p=5$ or $5\ |\ (p-1))$, ${\cal CGD}_{p^2}^2$ $(5\ |\
(p\pm 1))$ or ${\cal CGD}_{p^4}$ $(p=3$ or $p\geq 7)$. Let $A=\Aut(\G)$ and $N$ a nontrivial normal subgroup of $A$.

The graphs $K_6$, $FQ_4$ and ${\cal CGD}_{5^3}$ are basic because $N$ has one or two orbits by  MAGMA~\cite{BCP}. Since there is no pentavalent graph of odd order, by Proposition~\ref{QG}, $N$ is transitive on $V(\Gamma)$ or each orbit of $N$
on  $V(\Gamma)$ has odd length. In particular, ${\cal CD}_{p}$ $(p=5$ or $5\ |\ (p-1))$
are basic.

Let $\G={\cal CGD}_{p^2}^2$ $(5\ |\ (p\pm 1))$. It suffices to show that $|N|\not=p$. By MAGMA~\cite{BCP}, $\Aut({\cal
CGD}_{11^2}^2)$ has no normal subgroups of order $11$, and we may assume that $p>11$.
Suppose $|N|=p$. By Theorem~\ref{dip5-covers}, $A=\Aut({\cal
CGD}_{p^2}^2)=R(GD_{p^2})\rtimes D_5$, and hence $A/N$ has
stabilizer $D_5$ on $V(\G_N)$, which is impossible for $p>11$ by
Proposition~\ref{2p}. It follows that ${\cal
CGD}_{p^2}^2$ $(5\ |\ (p\pm 1))$ are basic. Finally let $\G={\cal
CGD}_{p^4}$ $(p=3$ or $p\geq 7)$. Similarly, $\Aut(\G)$ has no
normal subgroups of order $p$, $p^2$ or $p^3$ (for $p=11$ by
MAGMA~\cite{BCP} and for $p=3,7$ or $p\geq 13$ by
Theorem~\ref{dip5-covers} and Proposition~\ref{2p}), and hence ${\cal
CGD}_{p^4}$ $(p=3$ or $p\geq 7)$ are basic.

On the other hand, let $\G={\cal CGD}_{p^2}^1$ $(p=5$ or $5\ |\ (p-1))$, ${\cal
CGD}_{p^3}$ $(5\ |\ (p-1))$ or ${\cal CGD}_{p^4}$ $(p=2$ or $5)$
and let $A=\Aut(\G)$. To finish the proof, by Theorem~\ref{th-3.1} and Proposition~\ref{QG}, it
suffices to show that $A$ has a nontrivial normal subgroup having more than two
orbits. This is true for $\G={\cal CGD}_{5^2}^1, {\cal CGD}_{2^4}$
or ${\cal CGD}_{5^4}$ by MAGMA~\cite{BCP}, and moreover, $\G$ is a
normal cover of ${\cal CD}_5$, $FQ_4$ or ${\cal CGD}_{5^3}$,
respectively. Let $\G={\cal CGD}_{p^2}^1$ $(5\ |\ (p-1))$. Then
$A=\Aut({\cal CGD}_{p^2}^1)=R(GD_{p^2})\rtimes\mz_5$ and $A_u\cong\mz_5$ with $u\in V(\G)$.
The group $R(GD_{p^2})$ contains a characteristic subgroup of order $p^2$, say $\mz_p^2$, and each subgroup of $\mz_p^2$ is normal in
$R(GD_{p^2})$. The number of subgroups of order $p$ in $\mz_p^2$ is
$p+1$, and since $5\nmid (p+1)$, at least one of the subgroups of
order $p$ in $\mz_p^2$ is fixed by $A_u$. It follows that $A$ has a
normal subgroup of order $p$ and hence $\G$ is a normal cover of
${\cal CD}_{p}$ $(5\ |\ (p-1))$. Let $\G={\cal CGD}_{p^3}$ $(5\ |\
(p-1))$. Then $R(GD_{p^3})$ has a characteristic subgroup $\mz_p^3$, which has $p^2+p+1$ subgroups of order $p^2$. Since $5\nmid (p^2+p+1)$, $A$ has a normal
subgroup of order $p^2$ and $\G$ is a normal cover of ${\cal
CD}_{p}$ $(5\ |\ (p-1))$. This completes the proof. \hfill\qed

\section{ Pentavalent Symmetric graphs of order $2p^2$}

As an application of the results obtained in former sections, we classify pentavalent symmetric graphs of order $2p^2$.

\begin{theorem}\label{th-6.1}
Let $p$ be a prime and $\G$ a connected pentavalent symmetric graph
of order $2p^2$. Then $\G\cong {\cal CGD}_{p^2}^1$ $(p=5$ or $5\ |\
(p-1))$, ${\cal CGD}_{p^2}^2$ $(5\ |\ (p\pm 1))$ or ${\cal
CD}_{p^2}$ $(5\ |\ (p-1))$.
\end{theorem}

\demo Let $A=\Aut(\G)$ and $u\in V(\G)$. By McKay \cite{Mckay},
there is no pentavalent symmetric graph of order $8$ or $18$. Thus,
$p\geq 5$. By Proposition~\ref{stabilizer-val5}, $|A_u| \mid
2^9\cdot3^2\cdot 5$ and hence $|A|=2^{i}\cdot3^j\cdot5\cdot p^2$, $1\leq i\leq 10$, $0\leq j\leq 2$.

First we have an observation.\medskip

\f {\bf Observation:} If $p=5$ then $A$ has a non-abelian Sylow
$5$-group.\medskip

Suppose that $Q$ is an abelian Sylow $5$-subgroup of $A$. Since
$|V(\G)|=2\cdot 5^2$, we have $|Q|=5^3$ and $|Q_u|=5$. Thus, $Q$ has two
orbits, say $B_1$ and $B_2$ with $u\in B_1$. Since $Q$ is abelian,
$Q_u$ fixes every vertex in $B_1$ and all orbits of $Q_u$ in $B_2$
have length $5$. Thus, $\G$ contains a subgroup $K_{5,5}$, and by its connectedness, $\G\cong K_{5,5}$, which is contrary to 
$|\G|=2p^2$. The observation follows.\medskip

Let $H$ be a non-trivial normal abelian $p$-group of an
arc-transitive subgroup $B$ of $A$.
Then all orbits of $H$ have length $p$ or $p^2$ because $|\G|=2p^2$.
Suppose $H_u\not=1$. Since $H_u\lhd B_u$, we have $5\ |\ |H_u|$ and
$p=5$. By Proposition~\ref{QG}, $H$ has two orbits on $V(\G)$,
implying that $5^2\ |\ |H|$. Consequently, $5^3\ |\ |H|$ and $H$
is an abelian Sylow $5$-subgroup of $A$, which is impossible by
Observation. Thus, $H$ is semiregular on $V(\G)$, forcing that
$|H|=p$ or $p^2$. This implies that $H=\mz_p$, $\mz_{p^2}$ or
$\mz_p^2$. If $H=\mz_p^2$ then, by Theorem~\ref{dip5-covers}, $\G\cong
{\cal CGD}_{p^2}^1$ $(p=5$ or $5\ |\ (p-1))$ or ${\cal CGD}_{p^2}^2$
$(5\ |\ (p\pm 1))$, which are the first two families of graphs in the theorem. Thus, we may have the following assumption.
\medskip

\f {\bf Assumption:} Each non-trivial normal abelian $p$-group of
each arc-transitive subgroup of $A$ is semiregular on $V(\G)$ and
isomorphic to $\mz_p$ or $\mz_{p^2}$.\medskip

Let $N$ be a minimal normal subgroup of $A$. By Lemma~\ref{lem-3.2} and Assumption, $N\cong\mz_p$. Now we prove the following
claim.

\medskip
\f {\bf Claim:} $A$ has a semiregular subgroup $L$ of order $p^2$
such that $N\leq L$ and $\G$ is $N_A(L)$-arc-transitive.

By Proposition~\ref{QG}, the quotient graph $\G_N$ has order $2p$
and is $A/N$-arc-transitive. By Proposition~\ref{2p}, $\G_N\cong
{\cal CD}_{5}$ $( \cong K_{5,5})$ with $\Aut({\cal CD}_{5})\cong
(S_5\times S_5)\rtimes\mz_2$, $\G_N\cong {\cal CD}_{11}$ with
$\Aut({\cal CD}_{11})\cong \PGL(2, 11)$, or $\G_N\cong {\cal
CD}_{p}$ $(p>11$ and $5\ |\ (p-1))$ with $\Aut({\cal CD}_{p})\cong
D_p\rtimes\mz_5$. Since $A/N$ is arc-transitive, in each case $A/N$
contains a semiregular subgroup of order $p$, say $L/N$. Thus,
$10p\mid |A/N|$, $N\leq L$ and $L$ is a semiregular subgroup of
order $p^2$ in $A$.

For $p>11$ with $5\mid (p-1)$, we have $\Aut({\cal CD}_{p})\cong
D_p\rtimes\mz_5$. Since $10p\mid |A/N|$, we have $A/N=\Aut({\cal
CD}_{p})$ and hence $L$ is the unique normal Sylow $p$-subgroup of
$A$. Thus, $N_A(L)=A$ and $\G$ is $N_A(L)$-arc-transitive. For
$p=11$, $\Aut({\cal CD}_{11})=\PGL(2, 11)$ and $L/N$ is a Sylow $11$-subgroup of $\Aut({\cal CD}_{11})$. By ATLAS~\cite{atlas}, $\Aut({\cal CD}_{11})$ has a maximal subgroup $M/N\cong \mz_{11}\rtimes \mz_{10}$ such that $L/N\leq M/N$. Thus, $L\lhd M$. Since $M/N\not\leq \PSL(2,11)\leq \PGL(2,11)$, $M/N$ is vertex-transitive and hence arc-transitive on $\G_N$. It follows that $M$ is arc-transitive on $\G$ and $\G$ is $N_A(L)$-arc-transitive.
For $p=5$, $\G$ is a bipartite graph of order $5^2$ with the orbits of
$L$ as its bipartite sets. By Lemma~\ref{semiregulargraph}, $N_A(L)$
is vertex-transitive. Since Sylow $5$-subgroups of $A$ have order
$5^3$ and $|L|=p^2$, $N_A(L)$ contains a Sylow $5$-subgroup of $A$, and hence $\G$
is $N_A(L)$-arc-transitive.

Now we are ready to finish the proof. Set $B=N_A(L)$. By Claim, $\G$
is $B$-arc-transitive. Since $L\lhd B$ and $L$ has two orbits, $\G$
is bipartite with its bipartite sets as the orbits of $L$, say $B_1$
and $B_2$ with $u\in B_1$. By Assumption, $L\cong\mz_{p^2}$. Note
that $L$ is a Sylow $p$-subgroup of $C_B(L)$ for $p>5$ because
$p^3\nmid |A|$, where $C_B(L)$ is the centralizer of $L$ in $B$.
This is also true for $p=5$ by Observation. Thus, $L$ is a normal
Sylow $p$-subgroup of $C_B(L)$ and hence $C_B(L)=L\times K$, where
$K$ is a Hall $p'$-subgroup of $C_B(L)$. In particular, $K$ is
characteristic in $C_B(L)$. If $K\neq 1$, then $K\unlhd B$ because
$C_B(L)\unlhd B$, and since $p^2\nmid |K|$, $K$ has more than two
orbits on $V(\G)$. By Proposition~\ref{QG}, $K$ is semiregular and
hence $|K|=2$. But the quotient graph $\G_K$ would have odd
order $p^2$ and valency $5$, a contradiction. Thus, $K=1$ and
$C_B(L)=L$. Since $B/L=B/C_B(L)\lesssim \Aut(L)\cong \mz_{p(p-1)}$,
$B$ has a unique normal subgroup of order $2p^2$ containing $L$, say $R$.
By Lemma~\ref{semiregulargraph}, $B$ contains a regular dihedral
group of order $2p^2$ containing $L$. By the uniqueness of $R$, we have $R\cong D_{p^2}$ and hence
$\G$ is a Cayley graph on $R$. Since $R\lhd B$, Proposition~\ref{normal2p2} implies $\G\cong {\cal CD}_{p^2}$ ($5\ |\ (p-1)$). \hfill\qed

\medskip
\f {\bf Acknowledgement:} This work was supported by the National Natural Science Foundation of China (11571035, 11231008, 11271012) and by the 111 Project of China (B16002).

\end{document}